# Asymptotic Measures for Hyperbolic Piecewise Smooth Mappings of a Rectangle

Michael Jakobson* Sheldon Newhouse

April 9, 1996

**Abstract**

We prove the existence of Sinai-Ruelle-Bowen measures for a class of $C^2$ self-mappings of a rectangle with unbounded derivatives. The results can be regarded as a generalization of a well-known one dimensional Folklore Theorem on the existence of absolutely continuous invariant measures. In an earlier paper [8] analogous results were stated and the proofs were sketched for the case of invertible systems. Here we give complete proofs in the more general case of noninvertible systems, and, in particular, develop the theory of stable and unstable manifolds for maps with unbounded derivatives.

## 1 Folklore Theorem and SRB Measures

A well-known Folklore Theorem in one-dimensional dynamics can be formulated as follows.

**Folklore Theorem.** *Let $I = [0, 1]$ be the unit interval, and suppose $\{I_1, I_2, \ldots\}$ is a countable collection of disjoint open subintervals of $I$ such that $\bigcup_i I_i$ has the full Lebesgue measure in $I$. Suppose there are constants $K_0 > 1$ and $K_1 > 0$ and mappings $f_i : I_i \to I$ satisfying the following conditions.*

*partially supported by NSF Grant 9303369

*1. $f_i$ extends to a $C^2$ diffeomorphism from Closure$(I_i)$ onto $[0,1]$, and $\inf_{z\in I_i} \mid Df_i(z) \mid > K_0$ for all $i$.*

*2. $\sup_{z\in I_i} \dfrac{\mid D^2 f_i(z) \mid}{\mid Df_i(z) \mid} \mid I_i \mid < K_1$ for all $i$.*

*where $\mid I_i \mid$ denotes the length of $I_i$. Then, the mapping $F(z)$ defined by $F(z) = f_i(z)$ for $z \in I_i$, has a unique invariant ergodic probability measure $\mu$ equivalent to Lebesgue measure on $I$.*

For the proof of the Folklore theorem and the ergodic properties of $\mu$ see for example [2] and [14].

In an earlier paper [8] we presented an analog of this theorem for piecewise $C^2$ diffeomorphisms with unbounded derivatives with proof sketched. We now wish to give a more general version of the results in [8]. We refer the reader to that paper for relevant remarks and references.

Let $\tilde{Q}$ be a Borel subset of the unit square $Q$ in the plane $\mathbf{R}^2$ with positive Lebesgue measure, and let $F : \tilde{Q} \to \tilde{Q}$ be a Borel measurable map. An $F-$invariant Borel probablility measure $\mu$ on $Q$ is called a $Sinai - Ruelle - Bowen$ measure (or SRB-measure) for $F$ if $\mu$ is ergodic and there is a set $A \subset \tilde{Q}$ of positive Lebesgue measure such that for $x \in A$ and any continuous real-valued function $\phi : Q \to \mathbf{R}$, we have

$$\lim_{n\to\infty} \frac{1}{n} \sum_{k=0}^{n-1} \phi(F^k x) = \int \phi d\mu. \tag{1}$$

The set of all points $x$ for which (1) holds is called the *basin* of $\mu$.

Note that if $\mu$ is an SRB measure, and $m_1$ is the normalized Lebesgue measure on its basin, then the bounded convergence theorem gives the weak convergence of the averages $\frac{1}{n}\sum_{k=0}^{n-1} F_\star^k m_1$ of the iterates of $m_1$ to $\mu$. Hence, SRB measures occur as limiting mass distributions of sets of positive Lebesgue measure. This fact makes them natural objects to study.

We are interested in giving conditions under which certain two-dimensional maps $F$ which piecewise coincide with hyperbolic diffeomorphisms $f_i$ have SRB measures. As in the one-dimensional situation there is an essential difference between a finite and an infinite number of $f_i$. In the case of an infinite number of $f_i$, their derivatives grow with $i$ and relations between first and second derivatives become crucial.

## 2 Hyperbolicity and geometric conditions

Consider a countable collection $\xi = \{E_1, E_2, \dots, \}$ of full height closed curvilinear rectangles in $Q$. Assume that each $E_i$ lies inside a domain of definition of a $C^2$ diffeomorphism $f_i$ which maps $E_i$ onto its image $S_i \subset Q$. We assume each $E_i$ connects the top and the bottom of $Q$. Thus each $E_i$ is bounded from above and from below by two subintervals of the line segments $\{(x,y) : y = 1,\ 0 \le x \le 1\}$ and $\{(x,y) : y = 0,\ 0 \le x \le 1\}$. We assume that the left and right boundaries of $E_i$ are graphs of smooth functions $x^{(i)}(y)$ with $\left|\frac{dx^{(i)}}{dy}\right| \le \alpha$ where $\alpha$ is a real number satisfying $0 < \alpha < 1$. We further assume that the images $f_i(E_i) = S_i$ are narrow strips connecting the left and right sides of $Q$ and that they are bounded on the left and right by the two subintervals of the line segments $\{(x,y) : x = 0,\ 0 \le y \le 1\}$ and $\{(x,y) : x = 1,\ 0 \le y \le 1\}$ and above and below by the graphs of smooth functions $Y^i(X), \mid \frac{dY^{(i)}}{dX} \mid \le \alpha$. We will see later that the upper bounds on derivatives $\left|\frac{dx^{(i)}}{dy}\right| \le \alpha$ and $\left|\frac{dY^{(i)}}{dX}\right| \le \alpha$ follow from hyperbolicity conditions that we formulate below.

We call the $E_i's$ *posts* , the $S_i's$ *strips*, and we say the $E_i's$ are *full height* in $Q$ while the $S_i's$ are *full width* in $Q$.

For $z \in Q$, let $\ell_z$ be the horizontal line through $z$. We define $\delta_z(E_i) = diam(\ell_z \bigcap E_i), \delta_{i,max} = \max_{z\in Q} \delta_z(E_i), \delta_{i,min} = \min_{z\in Q} \delta_z(E_i)$.

We assume the following *geometric conditions*

G1 $int\ E_i \cap int\ E_j = \emptyset$ for $i \ne j$.

G2 $mes(Q \setminus \cup_i\ int\, E_i) = 0$ where $mes$ stands for Lebesgue measure,

G3 $-\sum_i \delta_{i,max} \log \delta_{i,min} < \infty$.

We emphasize that the strips $S_i$ can intersect in an arbitrary fashion, differently from condition G3 in ([8]).

In the standard coordinate system for a map $F : (x,y) \to (F_1(x,y), F_2(x,y))$ we use $DF(x,y)$ to denote the differential of $F$ at some point $(x,y)$ and $F_{jx}$, $F_{jy}$, $F_{jxx}$, $F_{jxy}$, etc., for partial derivatives of $F_j,\ j = 1,2$ .

Let $J_F(z) =\mid F_{1x}(z)F_{2y}(z) - F_{1y}(z)F_{2x}(z) \mid$ be the absolute value of the Jacobian determinant of $F$ at $z$.

*Hyperbolicity conditions.* There exist constants $0 < \alpha < 1$ and $K_0 > 1$ such that for each $i$ the map

$$F(z) = f_i(z) \text{ for } z \in E_i$$

satisfies

H1. $| F_{2x}(z) | + \alpha | F_{2y}(z) | + \alpha^2 | F_{1y}(z) | \leq \alpha | F_{1x}(z) |$

H2. $| F_{1x}(z) | - \alpha | F_{1y}(z) | \geq K_0$.

H3. $| F_{1y}(z) | + \alpha | F_{2y}(z) | + \alpha^2 | F_{2x}(z) | \leq \alpha | F_{1x}(z) |$

H4. $| F_{1x}(z) | - \alpha | F_{2x}(z) | \geq J_F(z) K_0$.

For a real number $0 < \alpha < 1$, we define the cones

$$K^u_\alpha = \{(v_1, v_2) : | v_2 | \leq \alpha | v_1 |\}$$

$$K^s_\alpha = \{(v_1, v_2) : | v_1 | \leq \alpha | v_2 |\}$$

and the corresponding cone fields $K^u_\alpha(z), K^s_\alpha(z)$ in the tangent spaces at points $z \in \mathbf{R}^2$.

Unless otherwise stated, we use the *max* norm on $\mathbf{R}^2$, $| (v_1, v_2) | = \max(| v_1 |, | v_2 |)$.

The following simple proposition relates conditions H1-H4 above with the usual definition of hyperbolicity in terms of cone conditions. It shows that conditions H1 and H2 imply that the $K^u_\alpha$ cone is mapped into itself by $DF$ and expanded by a factor no smaller than $K_0$ while H3 and H4 imply that the $K^s_\alpha$ cone is mapped into itself by $DF^{-1}$ and expanded by a factor no smaller than $K_0$.

**Proposition 2.1** *Under conditions H1-H4 above, we have*

$$DF(K^u_\alpha) \subseteq K^u_\alpha \tag{2}$$

$$v \in K^u_\alpha \Rightarrow | DFv | \geq K_0 | v | \tag{3}$$

$$DF^{-1}(K^s_\alpha) \subseteq K^s_\alpha \tag{4}$$

$$v \in K^s_\alpha \Rightarrow |\ DF^{-1}v\ | \geq K_0|\ v\ | \tag{5}$$

**Proof.** H1 implies (2):
Let $v = (v_1, v_2) \in K^u_\alpha$. Then, $|\ v\ | = |\ v_1\ |$ since $\alpha < 1$ and $|\ v_2\ | \leq \alpha|\ v_1\ |$.
Write $DF(v_1, v_2) = (F_{1x}v_1 + F_{1y}v_2, F_{2x}v_1 + F_{2y}v_2) = (u_1, u_2)$.
Then, using H1, we have

$$\begin{aligned}
|\ u_2\ | &= |\ F_{2x}v_1 + F_{2y}v_2\ | \\
&\leq |\ F_{2x}\ ||\ v_1\ | + |\ F_{2y}\ |\alpha|\ v_1\ | \\
&\leq |\ v_1\ |(|\ F_{2x}\ | + |\ F_{2y}\ |\alpha) \\
&\leq |\ v_1\ |(\alpha|\ F_{1x}\ | - |\ F_{1y}\ |\alpha^2) \\
&\leq \alpha|\ F_{1x}v_1 + F_{1y}v_2\ | \\
&= \alpha|\ u_1\ |
\end{aligned}$$

proving (2).
H2 implies (3):
Now, let $v = (v_1, v_2)$ be a unit vector in $K^u_\alpha$, so that $|\ v\ | = |\ v_1\ | = 1$ and $|\ v_2\ | \leq \alpha$.
Using H2 and the fact that $DF(v) \in K^u_\alpha$, we have

$$\begin{aligned}
|\ DF(v)\ | &= |\ u_1\ | \\
&= |\ F_{1x}v_1 + F_{1y}v_2\ | \\
&\geq |\ F_{1x}\ | - \alpha|\ F_{1y}\ | \\
&\geq K_0
\end{aligned}$$

which is (3).
The proofs that H3 and H4 imply (4) and (5) are similar using the fact that

$$DF^{-1} = \frac{1}{J_z}\begin{pmatrix} F_{2y} & -F_{1y} \\ -F_{2x} & F_{1x} \end{pmatrix}$$

This completes our proof of Proposition 2.1.

**Remark.** In ([8]) different hyperolicity conditions were assumed which implied the invariance of cones and uniform expansion with respect to the sum norm $| v | = | v_1 | + | v_2 |$ (see [3] and [7] for related hyperbolicity conditions). The methods here can be adapted to work under the assumptions of ([8]).

The map

$$F(z) = f_i(z) \text{ for } z \in int\ E_i$$

is defined almost everywhere on $Q$. Let $\tilde{Q}_0 = \bigcup_i int\ E_i$, and, define $\tilde{Q}_n, n > 0$, inductively by $\tilde{Q}_n = \tilde{Q}_0 \bigcap F^{-1}\tilde{Q}_{n-1}$. Let $\tilde{Q} = \bigcap_{n\geq 0} \tilde{Q}_n$ be the set of points whose forward orbits always stay in $\bigcup_i int\ E_i$. Then, $\tilde{Q}$ has full Lebesgue measure in $Q$, and $F$ maps $\tilde{Q}$ into itself.

The hyperbolicity conditions H1–H4 imply the estimates on the derivatives of the boundary curves of $E_i$ and $S_i$ which we described earlier. They also imply that any intersection $f_i E_i \bigcap E_j$ is full width in $E_j$. Further, $E_{ij} = E_i \bigcap f_i^{-1} E_j$ is a full height subpost of $E_i$ and $S_{ij} = f_j f_i E_{ij}$ is a full width substrip in $Q$.

Given a finite string $i_0 \dots i_{n-1}$, indexed by non-negative integers, we define inductively

$$E_{i_0 \dots i_{n-1}} = E_{i_0} \bigcap f_{i_0}^{-1} E_{i_1 i_2 \dots i_{n-1}}.$$

Then, each set $E_{i_0 \dots i_{n-1}}$ is a full height subpost of $E_{i_0}$.

Analogously, for a string $i_{-n+1} \dots i_0$ indexed by non-positive integers, we define.

$$S_{i_{-n+1} \dots i_0} = f_{i_0}(S_{i_{-n+1} \dots i_{-1}} \bigcap E_{i_0})$$

and get that $S_{i_{-n+1} \dots i_0}$ is a full width strip in $Q$. It is easy to see that $S_{i_{-n+1} \dots i_0} = (f_{i_0} \circ f_{i_{-1}} \circ \dots \circ f_{i_{-n+1}})(E_{i_{-n+1} \dots i_0})$ and that $f_{i_0}^{-1}(S_{i_{-n+1} \dots i_0})$ is a full-width strip in $E_{i_0}$.

For infinite strings, we have the following Proposition.

**Proposition 2.2** *Any $C^1$ map $F$ satisfying the above geometric conditions G1-G3 and hyperbolicity conditions H1–H4 has a "topological attractor"*

$$\Lambda = \bigcup_{\ldots i_{-n+1}\ldots i_{-1} i_0} \bigcap_{k \geq 0} S_{i_{-k}\ldots i_0}.$$

*The infinite intersections $\bigcap_{k=0}^{\infty} S_{i_{-k}\ldots i_0}$ define $C^1$ curves $\gamma = y(x)$, $|dy/dx| \leq \alpha$ which are the unstable manifolds for the points of the attractor. The infinite intersections $\bigcap_{k=0}^{\infty} E_{i_0\ldots i_k}$ define $C^1$ curves $x(y)$, $|dx/dy| \leq \alpha$ which are the stable manifolds for the points of the attractor.*

Proposition 2.2 is a well known fact in hyperbolic theory. For example it follows from Theorem 1 in [3]. See also [10].

**Remark 2.3** The distortion condition D1 and *distortion estimates* below imply that if our maps $f_i$ are $C^2$, then the unstable manifolds are actually $C^2$. Similar conditions on the inverses of $f_i$ imply that the stable manifolds are $C^2$. There are analogous conditions (see section 6) to guarantee that the invariant manifolds are $C^r$ for $r \geq 2$.

**Remark 2.4** The union of the stable manifolds contains the above set $\tilde{Q}$ which has full measure in $Q$. The trajectories of all points in $\tilde{Q}$ converge to $\Lambda$. That is the reason to call $\Lambda$ a topological attractor, although $F$ is not typically a well-defined mapping on all of $\Lambda$. However the convergence of Birkhoff averages to the unique SRB measure is a much stronger property. Condition D1 is natural in this context and may be necessary for the existence of the SRB measure. At present, we need to assume condition G3. This is used to prove absolute continuity of the stable foliation as in Section 10. It also implies that our SRB measure has finite entropy. We do not know if condition G3 is actually necessary for our results.

## 3 Distortion conditions and the main theorem

As we have a countable number of domains the derivatives of $f_i$ grow. We will need to formulate certain assumptions on the second derivatives. Unless otherwise stated, we will use the norm $\mid v \mid = \max(\mid v_1 \mid, \mid v_2 \mid)$ on vectors $v = (v_1, v_2)$, and the associated distance function $d((x,y),(x_1,y_1)) = \max(\mid x - x_1 \mid, \mid y - y_1 \mid)$.

As above, for a point $z \in Q$, let $l_z$ denote the horizontal line through $z$, and if $E \subseteq Q$, let $\delta_z(E)$ denote the diameter of the horizontal section $l_z \bigcap E$. We call $\delta_z(E)$ the $z - width$ of $E$.

In given coordinate systems we write $f_i(x,y) = (f_{i1}(x,y), f_{i2}(x,y))$. We use $f_{ijx}, f_{ijy}, f_{ijxx}, f_{ijxy}$, etc. for partial derivatives of $f_{ij}, j = 1, 2$.

We define

$$\mid D^2 f_i(z) \mid = \max_{j=1,2,(k,l)=(x,x),(x,y),(y,y)} \mid f_{ijkl}(z) \mid .$$

Next we formulate distortion conditions. These will be used to control the fluctuation of the derivatives of iterates of $F$ along vectors in $K^u_\alpha$ as in Lemma 7.1 and Proposition 8.1 below.

Suppose there is a constant $C_0 > 0$ such that the following *distortion condition* holds

D1 $\quad \sup_{z \in E_i, i \geq 1} \dfrac{\mid D^2 f_i(z) \mid}{\mid f_{i1x}(z) \mid} \delta_z(E_i) < C_0.$

**Theorem 3.1** *Let $F$ be a piecewise smooth mapping as above satisfying the geometric conditions G1–G3, the hyperbolicity conditions H1–H4, and the distortion condition D1.*

*Then, $F$ has an SRB measure $\mu$ whose basin has full Lebesgue measure in $Q$. Moreover, the natural extension of the system $(F, \mu)$ is measure-theoretically isomorphic to a Bernoulli shift, $F$ has finite entropy with respect to the measure $\mu$, and we have the formula*

$$h_\mu(F) = \lim_{n \to \infty} \frac{1}{n} \log \mid DF^n(z) \mid \tag{6}$$

*where the latter limit exists for Lebesgue almost all $z$ and is independent of such $z$.*

**Remark 3.2** Formula (6) says that the entropy can actually be computed by taking the logarithmic growth rate of the norms of $DF^n(z)$ for almost all $z$. It is actually true that if $v$ is any unit vector in the $K^u_\alpha$ cone in the tangent space to such a $z$, then

$$h_\mu(F) = \lim_{n\to\infty} \frac{1}{n} \log \mid DF^n(z)(v) \mid \tag{7}$$

This last expression can easily be implemented numerically.

**Remark 3.3** If we assume that the interiors of the strips $S_i$ are disjoint, then $(F, \mu)$ itself is isomorphic to a Bernoulli shift, and the entropy formula

$$h_\mu(F) = \int log|D^u F| d\mu$$

holds where $D^u F(z)$ is the norm of the derivative of $F$ in the unstable direction at $z$.

**Acknowledgement:** We wish to thank Francois Ledrappier and Dan Rudolph for useful conversations during the preparation of this paper.

# 4 Some estimates of partial derivatives

We will need to use the Mean Value Theorem for various partial derivatives of the mappings $f_i$ at points near the domain $E_i$. Since the $E_i$ are not necessarily convex subsets of $\mathbf{R}^2$, it will be useful to have our maps $f_i$ extended to neighborhoods $\mathcal{E}_i$ of $E_i$ which contain $\bigcup_{z \in E_i} B_{C\delta(z)}(z)$ where $C$ is a fixed positive constant and $B_{C\delta(z)}(z)$ denotes the ball about $z$ of radius $C\delta(z)$. Using the proof of the Whitney extension theorem in [1] it is possible to show that there is an extension $\tilde{f}_i$ of $f_i$ to such a neighborhood which satisfies the same properties H1-H4, D1, with possibly different constants. We will assume henceforth that our maps $f_i$ have such extensions.

We collect here some estimates which follow from our assumptions.

Let $f(x,y) = (f_1(x.y), f_2(x,y))$ be one of our maps $f_i$ on $E_i$.

**Lemma 4.1** *For $z \in E_i$, we have the estimates*

$$\frac{\mid f_{1y}(z) \mid}{\mid f_{1x}(z) \mid} \leq \alpha \tag{8}$$

$$\frac{\mid f_{2x}(z) \mid}{\mid f_{1x}(z) \mid} \leq \alpha \tag{9}$$

$$\frac{\mid f_{2y}(z) \mid}{\mid f_{1x}(z) \mid} \leq \frac{1}{K_0^2} + \alpha^2 \tag{10}$$

**Proof.**

We have

$$Df_z = \begin{pmatrix} f_{1x} & f_{1y} \\ f_{2x} & f_{2y} \end{pmatrix}$$

and

$$Df_{fz}^{-1} = \frac{1}{J_z} \begin{pmatrix} f_{2y} & -f_{1y} \\ -f_{2x} & f_{1x} \end{pmatrix}$$

where $J_z = f_{1x}f_{2y} - f_{2x}f_{1y}$.

Using $Df_z \begin{pmatrix} 1 \\ 0 \end{pmatrix} \in K_\alpha^u$ and $Df_{fz}^{-1} \begin{pmatrix} 0 \\ 1 \end{pmatrix} \in K_\alpha^s$ immediately gives

$$\frac{\mid f_{2x} \mid}{\mid f_{1x} \mid} \leq \alpha, \quad \frac{\mid f_{1y} \mid}{\mid f_{1x} \mid} \leq \alpha.$$

Now, we know that $\mid Df_{fz}^{-1}\begin{pmatrix} 0 \\ 1 \end{pmatrix} \mid \geq K_0$ in the max norm, so $\frac{1}{\mid J_z \mid} \max(\mid f_{1y} \mid, \mid f_{1x} \mid) \geq K_0$.

Hence, either $\mid J_z \mid K_0 \leq \mid f_{1y} \mid$ or $\mid J_z \mid K_0 \leq \mid f_{1x} \mid$.

The first case gives

$$(\mid f_{1x} f_{2y} \mid - \mid f_{1y} f_{2x} \mid) K_0 \leq \mid f_{1y} \mid$$

or

$$\begin{aligned} \frac{\mid f_{2y} \mid}{\mid f_{1x} \mid} &\leq \frac{\mid f_{1y} \mid}{K_0 \mid f_{1x} \mid^2} + \frac{\mid f_{1y} f_{2x} \mid}{\mid f_{1x} \mid^2} \\ &\leq \frac{\alpha}{K_0 \mid f_{1x} \mid} + \alpha^2 \\ &\leq \frac{1}{K_0^2} + \alpha^2. \end{aligned}$$

Analogously, in the second case,

$$(\mid f_{1x} f_{2y} \mid - \mid f_{1y} f_{2x} \mid) K_0 \leq \mid f_{1x} \mid$$

or

$$\begin{aligned} \frac{\mid f_{2y} \mid}{\mid f_{1x} \mid} &\leq \frac{1}{K_0 \mid f_{1x} \mid} + \alpha^2 \\ &\leq \frac{1}{K_0^2} + \alpha^2 \end{aligned}$$

Thus, in any case, we have

$$\frac{\mid f_{2y} \mid}{\mid f_{1x} \mid} \leq \frac{1}{K_0^2} + \alpha^2. \ QED$$

We have assumed that our maps $f_i$ have extensions to neighborhoods $\mathcal{E}_i$ of $E_i$ with the following properties.

The map $f_i$ takes $\mathcal{E}_i$ onto a set $\tilde{S}_i \subset \mathbf{R}^2$ such that

$$B_{C\delta_{z(E_i)}}(z) \subset \mathcal{E}_i \text{ for } z \in E_i \tag{11}$$

and

$$f_i \text{ satisfies H1–H4, D1 on } \mathcal{E}_i \tag{12}$$

Any $C^1$ curve $\gamma(t)$ such that $\gamma'(t) \in K^u_\alpha$ for all $t$ will be called a $K^u_\alpha$ curve. Similarly, a $K^s_\alpha$ curve is a $C^1$ curve $\gamma(t)$ for which $\gamma'(t) \in K^s_\alpha$ for all $t$. In this paper, all of our $K^u_\alpha$ curves will actually be of class $C^2$, and this will be assumed without further mention.

**Lemma 4.2** *Let $f_i$ have an extension to the neighborhood $\mathcal{E}_i$ as above. Then, there is a constant $C_1 > 0$ independent of $i$ such that if $z$ and $w$ lie on a $K^u_\alpha$ curve in $\mathcal{E}_i$, then*

$$\frac{|\ f_{i1x}(z)\ |}{|\ f_{i1x}(w)\ |} \le exp\left(C_1 \frac{|\ z-w\ |}{\delta_z(E_i)}\right).$$

**Proof.**

Write $f = f_i$.

Since $|\ Df_z \begin{pmatrix} 1 \\ 0 \end{pmatrix}\ | = \max(|\ f_{1x}(z)\ |, |\ f_{2x}(z)\ |) \ge K_0$ and $|\ f_{2x}(z)\ | \le \alpha |\ f_{1x}(z)\ |$, we know that

$$|\ f_{1x}(z)\ | = |\ Df_z \begin{pmatrix} 1 \\ 0 \end{pmatrix}\ | \ge K_0 > 1$$

so, for $w$ near $z$, both $f_{1x}(z)$ and $f_{1x}(w)$ have the same sign. We assume this sign is positive (replace $f$ by $-f$ otherwise).

Since $f$ extends to the neighborhood $\mathcal{E}_i$, and, for some constant $C > 0$, this last set contains the balls of radius $C\delta_z(E_i) > 0$ about points $z$ in $E_i$, the mean value theorem gives us that if $|\ z - w\ | \le C\delta_z(E_i)$, then there is a $\tau$ on the line segment joining $z$ and $w$ such that

$$|\ \log f_{1x}(z) - \log f_{1x}(w)\ | \le |\ \tfrac{f_{1xx}(\mathcal{T})}{f_{1x}(\mathcal{T})}\ ||\ z-w\ | + |\ \tfrac{f_{1xy}(\mathcal{T})}{f_{1x}(\mathcal{T})}\ ||\ z-w\ |\alpha$$

or

$$\frac{|\ f_{1x}(z)\ |}{|\ f_{1x}(w)\ |} \le exp\left((1+\alpha)C_0 \frac{|\ z-w\ |}{\delta_{\mathcal{T}}(E_i)}\right)$$

using the distortion estimate D1.

Let $z = (x_0, y_0)$, let $z_r = (x_r, y_r)$ be the point of intersection of the horizontal line $\ell_z$ with the right boundary curve of $E_i$, and let $z_\ell = (x_\ell, y_\ell)$ be the point of intersection of the horizontal line $\ell_z$ with the left boundary curve of $E_i$. Since $w$ lies on a $K^u_\alpha$ curve containing $z$, the line $\ell^0$ through $z$ and $\tau$ has equation $y - y_0 = \beta(x - x_0)$ for some $\beta$ with $\mid \beta \mid \leq \alpha$. Also, since the right boundary curve of $E_i$ through $z_r$ is a $K^s_\alpha - curve$, it is contained between the lines $\ell^-_r : x - x_r = -\alpha(y - y_r)$ and $\ell^+_r : x - x_r = \alpha(y - y_r)$. Similar statements hold for the left boundary curve of $E_i$ and the lines $\ell^-_\ell : x - x_\ell = -\alpha(y - y_\ell)$ and $\ell^+_\ell : x - x_\ell = \alpha(y - y_\ell)$. Using the intersections of the lines $\ell^0, \ell^\pm_r, \ell^\pm_\ell$, an elementary argument gives that

$$\frac{1}{1+\alpha^2} \leq \frac{1}{1+\mid \beta \mid \alpha} \leq \frac{\delta_\mathcal{T}(E_i)}{\delta_z(E_i)} \leq \frac{1}{1-\mid \beta \mid \alpha} \leq \frac{1}{1-\alpha^2}.$$

This gives the desired estimate for $z, w$ with $\mid z - w \mid \leq C\delta_z(E_i)$.

To get the general estimate of the Lemma, we simply find a sequence $z_0 = z, z_1, \ldots z_j = w$ with $z_k \in E_i$, $\mid z_k - z_{k+1} \mid < C\delta_z(E_i)$, each $z_k$ on the same $K^u_\alpha$ curve, and $j$ dependent only on $\alpha, C$, and $\delta_z(E_i)$. Using the estimate for each pair $z_i, z_{i+1}$ then easily gives us the general estimate to complete the proof of the Lemma. QED.

In some of our arguments below, it will simplify matters if we can take the constant $K_0$ in (3) and (5) to be large. The next lemma shows that this can be arranged by replacing $F$ by a fixed finite power $F^t$ with $t > 0$.

**Lemma 4.3** *Suppose the maps $f_i$ satisfy (2), (3), (4), (5), and D1 on the neighborhoods*

$$\bigcup_{z \in E_i} B_{C\delta_{z(E_i)}}(z),$$

*and let $t > 0$ be a positive integer.*

*Then there are positive constants $C_0 = C_0(t), C_2 = C_2(t)$ such that the maps $f_{i_{t-1}} \circ \ldots \circ f_{i_0}$ satisfy (2), (3), (4), and (5) with $K_0$ replaced by $K_0^t$ and D1 with $C_0$ replaced by $C_0(t)$ on the neighborhoods*

$$\bigcup_{z \in E_{i_0 \cdots i_{t-1}}} B_{C_2(t)\delta_{z(E_{i_0 \cdots i_{t-1}})}}(z)$$

**Proof.**

The proof is by induction on the number of elements in the composition. We assume that it holds for compositions of length $t$ and prove it for those of length $t+1$.

Let $B_{C(i_0 \dots i_t)}$ denote the set

$$\bigcup_{z \in E_{i_0 \dots i_t}} B_{C_2(t+1)\delta_z(E_{i_0 \dots i_t})}(z).$$

From Lemma 4.2, we can choose a constant $C_2(t+1) \in (0, C_2(t)) \subset (0,1)$ so that if $w \in B_{C(i_0 \dots i_t)}$, then $f_{i_0}(w) \in B_{C(i_1 \dots i_t)}$.

It is clear that the maps $f_{i_t} \circ \dots \circ f_{i_0}$ satisfy (2), (3), (4), and (5) with $K_0$ replaced by $K_0^{t+1}$, so we only need to be concerned with the statement regarding D1.

If $E$ is a subset of $Q$, $z \in E$, and $f(x,y) = (f_1(x,y), f_2(x,y))$, we set

$$\Theta_z(f, E) = \max_{i=1,2} \frac{\mid D^2 f_i(z) \mid}{\mid f_{1x}(z) \mid} \delta_z(E)$$

Let $f = f_{i_t} \circ \dots \circ f_{i_1}, g = f_{i_0}, h = f \circ g, E_f = E_{i_1 \dots i_t}, E_g = E_{i_0}, E_h = E_{i_0 \dots i_t}$, and, for $z \in E_h$, write $\Delta f = \delta_{gz}(E_f), \Delta g = \delta_z(E_g)$, and $\Delta h = \delta_z(E_h)$. Also, write $\Theta(f) = \Theta_{gz}(f, E_f), \Theta(g) = \Theta_z(g, E_g), \Theta(h) = \Theta_z(h, E_h)$.

Let us first estimate the quotient

$$\frac{\mid g_{1x}(w) \mid}{\mid g_{1x}(z) \mid}$$

for any $w \in \ell_z \cap \mathcal{E}_{i_0}$.

Note that $g_{1x}(w)$ and $g_{1x}(z)$ have the same sign. We assume it is positive. The argument when it is negative is similar.

Letting $C_1$ be the constant in Lemma 4.2, if $w, \bar{w} \in \ell_z \cap \mathcal{E}_{i_0}$, we have

$$\frac{\mid g_{1x}(w) \mid}{\mid g_{1x}(\bar{w} \mid} \leq exp(2C_1) \tag{13}$$

We can connect $w$ to $z$ in $\ell_z \cap \mathcal{E}_{i_0}$ by a chain of points $w = w_0, w_1, \dots w_k = z$ where $\mid w_i - w_{i+1} \mid < C_2(1)\delta_z(E_g)$, and $k \leq \frac{3}{C_2(1)}$.

Hence, putting $\zeta = exp\left(\frac{6C_1}{C_2(1)}\right)$, we have

$$\frac{| g_{1x}(w) |}{| g_{1x}(z) |} \leq \prod_{0\leq i<k} \frac{| g_{1x}(w_i) |}{| g_{1x}(w_{i+1} |} \leq \zeta \tag{14}$$

Interchanging $z, w$ in the above argument gives $\frac{| g_{1x}(w) |}{| g_{1x}(z) |} \geq \zeta^{-1}$. From these two inequalities we get, for any $w, \tau \in \ell_z \cap \mathcal{E}_{i_0}$,

$$\begin{aligned} \zeta^{-2} &\leq \frac{| g_{1x}(\tau) |}{| g_{1x}(z) |} \frac{| g_{1x}(z) |}{| g_{1x}(w) |} \\ &= \frac{| g_{1x}(\tau) |}{| g_{1x}(w) |} \\ &\leq \zeta^2 \end{aligned}$$

By the Chain Rule for partial derivatives we have the following formulas for $i = 1, 2$

$$h_{ix} = f_{ix} g_{1x} + f_{iy} g_{2x} \ ; \ h_{iy} = f_{ix} g_{1y} + f_{iy} g_{2y} \tag{15}$$

$$\begin{aligned} h_{ixx} &= f_{ixx} g_{1x}^2 + f_{ixy} g_{2x} g_{1x} + f_{iyx} g_{1x} g_{2x} + f_{iyy} g_{2x}^2 \\ &\quad + f_{ix} g_{1xx} + f_{iy} g_{2xx} \end{aligned} \tag{16}$$

$$\begin{aligned} h_{ixy} &= f_{ixx} g_{1y} g_{1x} + f_{ixy} g_{2y} g_{1x} + f_{iyx} g_{1y} g_{2x} + f_{iyy} g_{2y} g_{2x} \\ &\quad + f_{ix} g_{1xy} + f_{iy} g_{2xy} \end{aligned} \tag{17}$$

$$\begin{aligned} h_{iyy} &= f_{ixx} g_{1y}^2 + f_{ixy} g_{2y} g_{1y} + + f_{iyx} g_{1y} g_{2y} + f_{iyy} g_{2y}^2 \\ &\quad + f_{ix} g_{1yy} + f_{iy} g_{2yy} \end{aligned} \tag{18}$$

Let $w \in B_{C(i_0 \ldots i_t)}$. Except where otherwise mentioned, we compute the partial derivatives below at $w$.

From (15) and Lemma 4.1, we get

$$| h_{1x} | \geq | f_{1x} g_{1x} | (1 - \alpha^2) \tag{19}$$

From (16), we have, for $i = 1, 2$,

$$\left|\frac{h_{ixx}(w)}{h_{1x}(w)}\Delta h\right| \;\leq\; (1-\alpha^2)^{-1}\left[\Theta(f)|\,g_{1x}(w)\,|\frac{\Delta h}{\Delta f} + 2\Theta(f)|\,g_{2x}(w)\,|\frac{\Delta h}{\Delta f} + \Theta(f)|\,g_{2x}(w)\,|\frac{\Delta h}{\Delta f} \right.$$
$$\left. + 3\Theta(g)\frac{\Delta h}{\Delta g}\right]$$

Since the $g-$image of a horizontal line is a $K^u_\alpha$ curve and the boundaries of $E_f$ are $K^s_\alpha$ curves, we can use the mean value theorem and a simple geometric estimate to get a constant $C_3(\alpha) > 0$ such that

$$|\,g_{1x}(\tau)\,|\Delta h \leq C_3(\alpha)\Delta f$$

for some point $\tau$ in $\ell_z \cap E_h$.

Putting all these estimates together gives

$$\left|\frac{h_{ixx}(w)}{h_{1x}(w)}\Delta h\right| \leq C_4(\alpha)(\Theta(f)\zeta^2 + \Theta(g))$$

Similar estimates can be given for the quantities $\left|\frac{h_{ixy}(w)}{h_{1x}(w)}\Delta h\right|, \left|\frac{h_{iyy}(w)}{h_{1x}(w)}\Delta h\right|$.

Thus, we simply define $C_0(t+1)$ so that it is larger than $C_4(\alpha)(C_0(t)\zeta^2 + C_0)$ and we have proved Lemma 4.3.

## 5 Families of Fiber contractions

Fiber contraction maps were defined in [7] to provide a tool in the analysis of smoothness of stable and unstable manifolds. We collect here certain facts about parametrized families of fiber contraction maps and related concepts.

Let $(X, d_1), (Y, d_2)$ be complete metric spaces and give $X \times Y$ the metric

$$d((x, y), (x', y')) = \max(d_1(x, x'), d_2(y, y')).$$

Let $\pi_1 : X \times Y \to X, \pi_2 : X \times Y \to Y$ be the natural projections.

A pair of maps $(F, f)$ is called a fiber contraction on $X \times Y$ if the following properties hold.

1. $f : X \to X$ and $F : X \times Y \to X \times Y$ are continuous maps.
2. $\pi_1 F = f \pi_1$.
3. There is a constant $0 < K < 1$ such that for $x \in X, y, y' \in Y$, we have

$$d(F(x, y), F(x, y')) \leq K d_2(y, y').$$

We call $f$ the *base map* and $F$ the *total map* of the fiber contraction $(F, f)$.

Let $f$ be a continuous self-map of the complete metric space $X$. We say the a point $x_0 \in X$ is an attracting fixed point of $f$ if for every $x \in X$, the sequence of iterates $x, f(x), f^2(x), \ldots$ converges to $x_0$ as $n \to \infty$. Clearly if such an $x_0$ exists, it must be the unique fixed point of $f$.

Let $A$ be a topological space and consider a family $\{f_\lambda\}_{\lambda \in A}$ of self-maps of the complete metric space $X$. We say that the family is continuous if the map $(\lambda, x) \to f_\lambda(x)$ from $A \times X$ to $X$ is continuous.

A family $\{f_\lambda\}$ of self-maps of $X$ is called a *uniform family of contractions* if

1. there is a constant $0 < K < 1$ such that, for all $\lambda, x, x'$,

$$d(f_\lambda x, f_\lambda x') \leq K d(x, x').$$

2. the family $\{f_\lambda\}$ is continuous.

We say that a family $\{(F_\lambda, f_\lambda)\}$ of fiber contractions is a *uniform family of fiber contractions* if

1. the fiber Lipschitz constants are uniformly less than 1. That is, there is a constant $0 < K < 1$ such that for any $\lambda, x, y, y'$

$$d(F_\lambda(x,y), F_\lambda(x,y')) \leq K d_2(y,y')$$

2. the families $\{F_\lambda\}$ and $\{f_\lambda\}$ are continuous.

The following Proposition is standard (see e.g. [6]) and its proof will be omitted.

**Proposition 5.1** *If $\{f_\lambda\}$ is a uniform family of contractions of the complete metric space $X$, and $x_\lambda$ is the fixed point of $f_\lambda$, then the family $\{x_\lambda\}$ depends continuously on $\lambda$.*

**Proposition 5.2** *Suppose $\{(F_\lambda, f_\lambda)\}$ is a uniform family of fiber contractions whose base maps $\{f_\lambda\}$ have attracting fixed points $\{x_\lambda\}$ depending continuously on $\lambda$. Then, each of the maps $F_\lambda$ has an attracting fixed point of the form $(x_\lambda, y_\lambda) \in X \times Y$ and the family $\{(x_\lambda, y_\lambda)\}$ depends continuously on $\lambda$.*

**Proof.** Letting $x_\lambda$ be the fixed point of the base map $f_\lambda$, Hirsch and Pugh prove in [7] that $F_\lambda$ has an attracting fixed point of the form $(x_\lambda, y_\lambda)$ where $y_\lambda$ is the fixed point of the map $F(x_\lambda, \cdot)$ on $Y$. Since $x_\lambda$ depends continuously on $\lambda$, the family $\{F(x_\lambda, \cdot)\}$ is family of uniform contractions on $Y$. Therefore, by Proposition 5.1, the fixed points $\{y_\lambda\}$ depend continuously on $\lambda$. QED.

The following corollary is proved by induction using Propositions 5.1 and 5.2.

**Corollary 5.3** *Suppose $X_1 \times X_2 \times \ldots \times X_N$ is a sequence of complete metric spaces and $\{F_{\lambda,i}\}, 1 \leq i \leq N$ is a sequence of maps with the following properties.*

1. *$\{F_{\lambda,1}\}$ is a uniform family of contractions on $X_1$.*

*2. For* $2 \leq i \leq N$, $\{(F_{\lambda,i}, F_{\lambda,i-1}\}$ *is a uniform family of fiber contractions on* $\prod_{1\leq j\leq i} X_j$.

*Then, each of the families* $\{F_{\lambda,i}\}$ *has an attracting family of fixed points* $\{x_{\lambda,i}\}$ *which depends continuously on* $\lambda$.

## 6 Invariant Manifolds

We consider the collection $\xi = \{E_1, E_2, \ldots\}$ of rectangles as above and the sequence $(f_1, f_2, \ldots)$ of $C^2$ diffeomorphisms with $f_i(E_i) = S_i$ satisfying G1–G3, H1–H4, and D1. From Proposition 2.1, using the max norm on $\mathbf{R}^2$, we have, for each $i$,

$$Df_i(K^u_\alpha) \subseteq K^u_\alpha \tag{20}$$

$$v \in K^u_\alpha \Rightarrow |\, Df_i v \,| \geq K_0 |\, v \,| \tag{21}$$

$$Df_i^{-1}(K^s_\alpha) \subseteq K^s_\alpha \tag{22}$$

$$v \in K^s_\alpha \Rightarrow |\, Df_i^{-1} v \,| \geq K_0 |\, v \,| \tag{23}$$

For each finite sequence $i_{-n+1} \ldots i_0 \ldots i_{n-1}$ we have defined, in Section 2, the sets $E_{i_0 \ldots i_{n-1}}, S_{i_{-n+1} \ldots i_0}$.

Given a non-positive itinerary $\mathbf{i} = (\ldots i_{-n} i_{-n+1} \ldots i_0)$, we consider the set $W^u_{\mathbf{i}} = E_{i_0} \cap \bigcap_{n \geq 0} S_{i_{-n} \ldots i_{-1}}$. Clearly, $W^u_{\mathbf{i}}$ is a closed, connected full-width subset of $E_{i_0}$. Its image $FW^u_{\mathbf{i}} = f_{i_0} W^u_{\mathbf{i}}$ is the set $\bigcap_{n \geq 0} S_{i_{-n} \ldots i_0}$, a full-width connected subset of $Q$. The next result shows that $FW^u_{\mathbf{i}}$ is a $C^2$ curve which depends continuously on $\mathbf{i}$.

For convenience, we let $D^0 \psi = \psi$ for a function $\psi$.

**Theorem 6.1** *There is a constant $K > 0$ such that for each non-positive itinerary $\mathbf{i} = (\ldots i_{-n} \ldots i_0)$, the set $FW^u_{\mathbf{i}}$ is the graph of a $C^2$ function $g_{\mathbf{i}} : I \to I$ such that, for $z \in I$,*

$$|\, Dg_{\mathbf{i}}(z) \,| \leq \alpha \tag{24}$$

*and*

$$|\, D^2 g_{\mathbf{i}}(z) \,| \leq K. \tag{25}$$

*Further, given $\epsilon > 0$, there is a positive integer $N > 0$ such that if $\mathbf{i} = (\ldots i_{-n} \ldots i_0)$ and $\mathbf{j} = (\ldots j_{-n} \ldots j_0)$ are non-positive itineraries with $i_{-\ell} = j_{-\ell}$ for $0 \leq \ell \leq N$, then*

$$| D^k g_{\mathbf{i}}(z) - D^k g_{\mathbf{j}}(z) | < \epsilon \tag{26}$$

*for* $z \in I$ *and* $0 \le k \le 2$.

**Remark 6.2** The proof of Theorem 6.1 uses graph transform techniques as in [7],[12]. However, since our maps have unbounded derivatives, and the off-diagonal terms of our derivatives are not small, certain modifications of the techniques in [7], [12] are necessary.

It can be shown that if $f_i$ is $C^r$ for $r \ge 2$, then the curves $W^u_{\mathbf{i}}$ are $C^r$ and depend continuously on $\mathbf{i}$ in the $C^r$ sense provided the $f_i$ satisfy the $r - th$ order distortion condition

$$\sup_{z \in \mathcal{E}_i, i \ge 1, 2 \le k \le r} \frac{| D^k f_i(z) |}{| f_{i1x}(z) |^{k-1}} \delta_z(\mathcal{E}_i) < C_0$$

where $D^k f_i(z)$ is the supremum of the $\ell - th$ order partial derivatives of $f_i$ at $z$ for $\ell \le k$

We proceed toward the proof of Theorem 6.1.

Notice that if we replace $F = \{f_i\}$ by a positive power $F^t, t > 0$, and $\xi$ by the collection $\{E_{i_0 \ldots i_{t-1}}\}$, we may assume that $K_0$ is as large as we wish in (21),(23). In the present section we will take $K_0 > 4$. Of course this changes the distortion constant $C_0$ in D1 to some $C_1 = C_1(f, t)$ but this will not cause us difficulties.

Let $\mathbf{N}$ be the set of positive integers, and let $\Sigma = \mathbf{N}^{\mathbf{Z}}$ be the space of doubly infinite sequences $\mathbf{i} = (\ldots i_{-1} i_0 i_1 \ldots)$ of elements of $\mathbf{N}$ with the product topology. Let $\sigma : \Sigma \to \Sigma$ be the usual left shift automorphism.

For an element $\mathbf{i} \in \Sigma$, let $\mathbf{i}^+ = (i_0 i_1 \ldots)$ be its non-negative part, and let $\mathbf{i}^- = (\ldots i_{-1} i_0)$ be its non-positive part. Set $W^s_{\mathbf{i}^+} = \bigcap_{n \ge 0} E_{i_0 \ldots i_n}$ and $W^u_{\mathbf{i}^-} = E_{i_0} \cap \bigcap_{n \ge 0} S_{i_{-n} \ldots i_{-1}}$.

It follows from (20)–(23) that the sets $W^s_{\mathbf{i}^+}, W^u_{\mathbf{i}^-}$ intersect in a unique point and there is a continuous map $\pi : \Sigma \to Q$ defined by

$$\{\pi((\ldots i_{-1} i_0 i_1 \ldots))\} = W^s_{\mathbf{i}^+} \bigcap W^u_{\mathbf{i}^-}$$

Moreover, for each $\mathbf{i} \in \Sigma$ there is a splitting $T_{\pi(\mathbf{i})} \mathbf{R}^2 = E^u_{\pi(\mathbf{i})} \oplus E^s_{\pi(\mathbf{i})}$ which depends continuously on $\mathbf{i}$ and is such that $Df_{i_0}$ maps $E^u_{\pi(\mathbf{i})}$ to $E^u_{\pi(\sigma \mathbf{i})}$ and $E^s_{\pi(\mathbf{i})}$ to $E^s_{\pi(\sigma\mathbf{i})}$. The arguments for these facts are analogous to standard arguments in hyperbolic theory (e.g., to prove that $C^1$ perturbations of the Smale horseshoe diffeomorphism have a hyperbolic non-wandering set) and will not be given here.

Thus, the matrix of $DF$ is diagonal with respect to the splitting $E^u \oplus E^s$ on the image of $\pi$.

For $z = \pi(\mathbf{i})$ and $v \in T_z\mathbf{R}^2$, we write $v = (v_1, v_2) \in E^u_z \oplus E^s_z$ and define $\mid v \mid = \mid v \mid_z = \max(\mid v_1 \mid, \mid v_2 \mid)$. This norm depends continuously on $\mathbf{i} \in \Sigma$.

We will identify all tangent spaces with the space $\mathbf{R}^2$ itself by standard translations.

It will be convenient to use the subundles $E^u, E^s$ to define affine local coordinates near points $z, f_i z$ in which $Df_{iz}$ becomes diagonal and in which $Df_{iw}$ is nearly diagonal for $\mid w - z \mid$ no larger than a fixed multiple of $\delta_z(E_i)$. Here $i = i_0$ with $z = \pi(\mathbf{i})$.

Toward this end, let $A_z$ be the affine automorphism of $\mathbf{R}^2$ such that

1. $A_z(z) = z$.

2. $DA_z \begin{pmatrix} 1 \\ 0 \end{pmatrix} = \begin{pmatrix} 1 \\ a_z \end{pmatrix} \in E^u_z$.

3. $DA_z \begin{pmatrix} 0 \\ 1 \end{pmatrix} = \begin{pmatrix} b_z \\ 1 \end{pmatrix} \in E^s_z$.

Since $E^u_z \subseteq K^u_\alpha$ and $E^s_z \subseteq K^s_\alpha$, we have $\mid a_z \mid \le \alpha, \mid b_z \mid \le \alpha$.

Let $\tilde{f} = \tilde{f}_i = A^{-1}_{f_i z} f_i A_z$ be the local representative of $f_i$ using the $A_{f_i z}, A_z$ coordinates. Note that $\tilde{f}$ is defined on the affine image $A^{-1}_z(\mathcal{E}_i)$ of $\mathcal{E}_i$.

Then the matrix $D\tilde{f}_z$ is diagonal. For $w$ near $z$ in $A_z(\mathcal{E}_i)$, let

$$D\tilde{f}_w = \begin{pmatrix} \tilde{f}_{1x}(w) & \tilde{f}_{1y}(w) \\ \tilde{f}_{2x}(w) & \tilde{f}_{2y}(w) \end{pmatrix}$$

and set

$$\epsilon_{12}(w) = \frac{\mid \tilde{f}_{1y}(w) \mid}{\mid \tilde{f}_{1x}(w) \mid},\ \epsilon_{21}(w) = \frac{\mid \tilde{f}_{2x}(w) \mid}{\mid \tilde{f}_{1x}(w) \mid},\ \epsilon_{22}(w) = \frac{\mid \tilde{f}_{2y}(w) \mid}{\mid \tilde{f}_{1x}(w) \mid}.$$

We wish to estimate $\epsilon_{ij}(w)$ for $w$ near $z$ in $\mathcal{E}_i$. It follows from the definitions that $\epsilon_{ij}(z) = 0, i \ne j$. Also, (21) and (23) imply $\epsilon_{22}(z) \le \frac{1}{K_0^2} < \frac{1}{16}$ since $K_0 > 4$.

**Lemma 6.3** *There are constants* $C_2 \in (0,1), C_3 > 0, C_4 > 0$, *such that for* $z \in E_i, w \in A_z^{-1}(\mathcal{E}_i)$, *if* $\mid w - z \mid < C_2\delta_z(E_i)$, *then*

$$\mid \epsilon_{ij}(w) - \epsilon_{ij}(z) \mid \leq C_3 \frac{\mid z - w \mid}{\delta_z(E_i)}, \tag{27}$$

*and*

$$C_4 \frac{1}{\delta_z(E_i)} \leq \mid \tilde{f}_{1x}(w) \mid \leq C_4^{-1} \frac{1}{\delta_z(E_i)} \tag{28}$$

**Proof.**

To begin with, let us choose $C_2 \in (0,1)$ so that if $\mid w - z \mid < C_2\delta_z(E_i)$, then $w \in \mathcal{E}_i \cap A_z^{-1}(\mathcal{E}_i)$ and $\tilde{f}$ satisfies D1 for some (possibly different) constant $C_0$. Since $A_{f_i z}^{-1}$ and $A_z$ are uniformly bounded, it is possible to choose $C_0$ and $C_2$ independent of $z \in E_i$ and $i \geq 1$.

We next show that there are constants $C_5 > 0, C_6 > 0$ such that for $z \in E_i$ and $\mid z - w \mid < C_5\delta_z(E_i)$,

$$C_6^{-1} \leq \frac{\mid f_{1x}(w) \mid}{\mid \tilde{f}_{1x}(w) \mid} \leq C_6. \tag{29}$$

Since

$$\begin{aligned} D\tilde{f}_w \begin{pmatrix} 1 \\ 0 \end{pmatrix} &= \begin{pmatrix} \tilde{f}_{1x}(w) \\ \tilde{f}_{2x}(w) \end{pmatrix} \\ &= DA_{fz}^{-1} Df_{A_z w} DA_z \begin{pmatrix} 1 \\ 0 \end{pmatrix} \\ &= \frac{1}{J_{fz}} \begin{pmatrix} 1 & -b_{fz} \\ -a_{fz} & 1 \end{pmatrix} \begin{pmatrix} f_{1x} & f_{1y} \\ f_{2x} & f_{2y} \end{pmatrix} \begin{pmatrix} 1 & b_z \\ a_z & 1 \end{pmatrix} \end{aligned}$$

where $J_{fz} = 1 - a_{fz}b_{fz}$ and the partial derivatives of $f_1, f_2$ are evaluated at $A_z(w)$, we have

$$\begin{aligned} J_{fz}\tilde{f}_{1x}(w) &= f_{1x} + f_{1y}a_z - b_{fz}f_{2x} - b_{fz}f_{2y}a_z \\ J_{fz}\tilde{f}_{1y}(w) &= f_{1x}b_z + f_{1y} - b_{fz}f_{2x}b_z - b_{fz}f_{2y} \\ J_{fz}\tilde{f}_{2x}(w) &= -a_{fz}f_{1x} - a_{fz}f_{1y}a_z + f_{2x} + f_{2y}a_z \end{aligned}$$

$$J_{fz}\tilde{f}_{2y}(w) \;=\; -a_{fz}f_{1x}b_z - a_{fz}f_{1y} + f_{2x}b_z + f_{2y}$$

Using the first equation above, the fact that $\mid J_{fz} \mid \geq 1-\alpha^2$, and the estimates (8), (9), (10) at $A_z(w)$ we get

$$\mid \tilde{f}_{1x}(w) \mid \leq C \mid f_{1x}(A_z(w)) \mid$$

for some constant $C$. But, from Lemma 4.2 we have $\mid f_{1x}(A_z(w)) \mid$ is bounded above by const$\mid f_{1x}(z) \mid$, so this gives the lower bound in (29).

For the upper bound, we will obtain the two estimates

$$\mid \tilde{f}_{1x}(w) \mid \geq C \mid D\tilde{f}_w \begin{pmatrix} 1 \\ 0 \end{pmatrix} \mid \tag{30}$$

and

$$\mid D\tilde{f}_w \begin{pmatrix} 1 \\ 0 \end{pmatrix} \mid \geq C \mid f_{1x}(w) \mid \tag{31}$$

for some constant $C > 0$.

To prove (30), we note first note that the vector $v = \begin{pmatrix} 1 & b_z \\ a_z & 1 \end{pmatrix} \begin{pmatrix} 1 \\ 0 \end{pmatrix}$ is in the cone $K^u_\alpha$. Since $Df_{A_z}$ preserves this cone, we have that $Df_{A_z}(v)$ is a constant multiple of the vector $\begin{pmatrix} 1 \\ \bar{a} \end{pmatrix}$ for some $\bar{a}$ with $\mid \bar{a} \mid \leq \alpha$.

Thus,

$$D\tilde{f}_w \begin{pmatrix} 1 \\ 0 \end{pmatrix} = \begin{pmatrix} \tilde{f}_{1x}(w) \\ \tilde{f}_{2x}(w) \end{pmatrix}$$

is a constant multiple of the vector

$$\begin{pmatrix} 1 & -b_{fz} \\ -a_{fz} & 1 \end{pmatrix} \begin{pmatrix} 1 \\ \bar{a} \end{pmatrix} = \begin{pmatrix} 1 - b_{fz}\bar{a} \\ -a_{fz} + \bar{a} \end{pmatrix}$$

This gives

$$\frac{\mid \tilde{f}_{2x}(w) \mid}{\mid \tilde{f}_{1x}(w) \mid} \leq \frac{2\alpha}{1-\alpha^2}$$

and, hence,

$$
\begin{aligned}
\mid \begin{pmatrix} \tilde{f}_{1x}(w) \\ \tilde{f}_{2x}(w) \end{pmatrix} \mid &= \max(\mid \tilde{f}_{1x}(w) \mid, \mid \tilde{f}_{2x}(w) \mid) \\
&\leq \max(1, \frac{2\alpha}{1-\alpha^2}) \mid \tilde{f}_{1x}(w) \mid
\end{aligned}
$$

and (30) follows.

Next we go to the proof of (31).

Since the matrix $\begin{pmatrix} 1 & b_{fz} \\ a_{fz} & 1 \end{pmatrix}$ and its inverse are uniformly bounded as are $J_{fz}$ and its inverse, we have

$$
\mid D\tilde{f}_w \begin{pmatrix} 1 \\ 0 \end{pmatrix} \mid \geq C \mid Df_{A_z(w)} \begin{pmatrix} 1 \\ a_z \end{pmatrix} \mid
$$

But,

$$
Df_{A_z(w)} \begin{pmatrix} 1 \\ a_z \end{pmatrix} = \begin{pmatrix} A_1 \\ A_2 \end{pmatrix}
$$

where $A_1 = f_{1x} + a_z f_{1y}$.

so,

$$
\begin{aligned}
\mid Df_{A_z(w)} \begin{pmatrix} 1 \\ a_z \end{pmatrix} \mid &\geq \mid f_{1x} \mid - \alpha^2 \mid f_{1x} \mid \\
&\geq (1-\alpha^2) \mid f_{1x} \mid \\
&\geq (1-\alpha^2) \mid f_{1x}(A_z^{-1}w) \mid \\
&\geq C \mid f_{1x}(w) \mid
\end{aligned}
$$

This completes the proof of (29).

Next, we give the proof of the estimate

$$
\mid \epsilon_{12}(w) - \epsilon_{12}(z) \mid \leq C_3 \frac{\mid z-w \mid}{\delta_z(E_i)} \tag{32}
$$

The other estimates for (27) are similar.

Since $\tilde{f}_{1y}(z) = 0$, we need to estimate $\mid \frac{\tilde{f}_{1y}(w)}{\tilde{f}_{1x}(w)} \mid$.

But,

$$J_{fz}\tilde{f}_{1y}(w) = f_{1x}(A_z w)b_z + f_{1y}(A_z w) - b_{fz}f_{2x}(A_z w)b_z - b_{fz}f_{2y}(A_z w),$$

so,

$$\mid \tilde{f}_{1y}(w) \mid \leq C \max_{i,j,k,\tau} \mid f_{ijk}(\tau) \mid\mid A_z(w) - z \mid.$$

Now, we know that the quantities $\frac{\delta_{z(E_i)}}{\delta_{\tau(E_i)}}, \frac{\mid f_{1x}(w) \mid}{\mid \tilde{f}_{1x}(w) \mid}$ are bounded above and below, and, by Lemma (4.2), the same holds for $\frac{\mid f_{1x}(\tau) \mid}{\mid \tilde{f}_{1x}(w) \mid}$. This gives (32) and (27).

For (28), notice that $f(\ell_z \bigcap E_i)$ is a full-width $K^u_\alpha$ curve in $Q$. In the max metric, it has unit length. By the Mean Value Theorem there is a $\tau \in \ell_z \bigcap E_i$ such that

$$\mid Df(\tau) \begin{pmatrix} 1 \\ 0 \end{pmatrix} \mid \delta_z(E_i) = 1$$

But,

$$\begin{aligned} \mid Df(\tau) \begin{pmatrix} 1 \\ 0 \end{pmatrix} \mid &= \max(\mid f_{1x}(\tau) \mid, \mid f_{2x}(\tau) \mid) \\ &= \mid f_{1x}(\tau) \mid \end{aligned}$$

so,

$$\mid f_{1x}(\tau) \mid = \frac{1}{\delta_z(E_i)}.$$

Since, $\frac{\mid f_{1x}(\tau) \mid}{\mid \tilde{f}_{1x}(w) \mid}$ is bounded above and below, (28) follows.

This completes the proof of Lemma 6.3.

For $\epsilon > 0$, let $B_\epsilon(z) = \{w \in \mathbf{R}^2 : \mid w - z \mid_z \leq \epsilon\}$. Here $\mid w - z \mid_z$ refers to the max norm in the image of the affine coordinate map $A_z$. The set $B_\epsilon(z)$ is then a parallelogram centered at $z$ with sides parallel to $E^u_z, E^s_z$. Write $B_\epsilon(z) = B^u_\epsilon(z) \times B^s_\epsilon(z)$ where $B^u_\epsilon(z)$ is a line segment centered at $z$ parallel to $E^u_z$, and $B^s_\epsilon(z)$ is a line segment centered at $z$ parallel to $E^s_z$.

A *full-width curve of slope less than 1* in $B_\epsilon(z)$ is the graph of a function $\phi : B^u_\epsilon(z) \to B^s_\epsilon(z)$ in which $\phi$ is Lipschitz with Lipschitz constant less than 1.

With $z = \pi(\mathbf{i})$, let $z_0 = z, z_j = \pi(\sigma^j \mathbf{i})$ for $j \leq 0$.

Our next goal, as is usual in invariant manifold theory, is to find a sequence of numbers $\epsilon_j > 0$ such that the neighborhoods $B_j = B_{\epsilon_j}(z_j)$ have the following properties.

B1 If $z_j \in E_{i_j}$, then $B_j \subseteq \mathcal{E}_{i_j}$.

B2 $\tilde{f}_{i_j}(B_j)$ overflows $B_{j+1}$ in the sense that if $\gamma$ is a full-width curve of slope less than 1 in $B_j$ passing through $z_j$ , then $\tilde{f}_{i_j}(\gamma) \bigcap B_{j+1}$ is a full-width curve of slope less than 1 in $B_{j+1}$ passing through $z_{j+1}$.

Let $\bar{\epsilon}_j = C_2 \delta_{z_j}(E_{i_j})$ where $C_2$ is the constant of Lemma 6.3. By Lemma 6.3, for $| w - z_j | < \bar{\epsilon}_j$, the matrix of $D\tilde{f}_{i_j}(w)$ is hyperbolic with off-diagonal terms small compared to $| \tilde{f}_{1x}(w) |$. This implies that the image $\tilde{f}_{i_j}\gamma$ of a curve $\gamma$ as above will have slope less than 1 in $B_{\bar{\epsilon}_{j+1}}(z_{j+1})$. Letting $f = \tilde{f}_{i_j}$, and using $a \sim b$ to mean $\frac{a}{b}$ is bounded above and below, we have $length(f\gamma) \sim | \tilde{f}_{1x}(\tau) | C_2 \delta_{z_j}(E_{i_j})$ for some $\tau \in \gamma$. In the proof of Lemma 6.3 we saw that $\delta_{z_j}(E_{i_j}) \sim \frac{1}{| f_{1x}(\tau_1) |} \sim \frac{1}{| \tilde{f}_{1x}(\tau_1) |}$ and $| \tilde{f}_{1x}(\tau) | \sim | \tilde{f}_{1x}(z_j) | \sim | \tilde{f}_{1x}(\tau_1) |$. It follows that $\tilde{f}_{i_j}\gamma$ contains a neighborhood of fixed size $C_7$ about $\tilde{f}_{i_j}(z_j)$ in $\tilde{f}_{i_j}\gamma$.

Let

$$\epsilon_j = \begin{cases} \bar{\epsilon}_j \text{ if } \bar{\epsilon}_j < C_7 \\ C_7 \text{ if } \bar{\epsilon}_j > C_7 \end{cases}$$

Then, the overflowing property above is satisfied.

Now fix a non-positive itinerary $\mathbf{i} = (\ldots i_{-n} \ldots i_0)$. We first show that $W^u_{\mathbf{i}}$ contains the graph of a $C^2$ function $g_{\mathbf{i}} : B^u_0 \to B^s_0$ such that, for all $w \in B^u_0$

$$| Dg_{\mathbf{i}}(w) | \leq 1 \tag{33}$$

and

$$| D^2 g_{\mathbf{i}}(w) | \leq K_2 \tag{34}$$

where $K_2$ is independent of $\mathbf{i}$.

We also will show that the functions $g_{\mathbf{i}}$ depend continuously on $\mathbf{i}$.

Once these things are done, the proof of Theorem 6.1 is completed as follows.

Let $\mathbf{j} = (\ldots i_{-1} i_0 j_1 j_2 \ldots)$ be a doubly infinite itinerary which agrees with $\mathbf{i}$ for non-positive indices. Let $z_0 = \pi(\mathbf{j})$. Then there is a $k > 0$ independent of $\mathbf{i}, \mathbf{j}$ such that $f_{i_{-k}}^{-1} \circ \ldots \circ f_{i_{-1}}^{-1}(W_{\mathbf{i}}^u) \subseteq B_{-k}$. Note that here we use the original maps $f_{i_j}$, not the affine representatives $\tilde{f}_{i_j}$.

Thus $W_{\mathbf{i}}^u$ is the $f_{i_{-1}} \circ \ldots \circ f_{i_{-k}}$−image of a curve of bounded slope and bounded $C^2$ size. Letting $F^k = f_{i_{-1}} \circ \ldots \circ f_{i_{-k}}$ we have that $W_{\mathbf{i}}^u$ is the graph of a function $\Gamma(F^k, g)$ where $F^k$ has bounded distortion and $g$ has bounded $C^1, C^2$ sizes. Using the formulas (36), (37), and (38) which appear in the second derivative of the graph transform function then gives that $\Gamma(F^k, g)$ also has bounded $C^2$ size. The same argument then works for $\Gamma(F^{k+1}, g)$ and this gives (25). A similar argument gives the continuity statement in Theorem 6.1.

To get estimate (24) first note that hyperbolicity conditions imply that any vector $v$ in the tangent space to a point in $W_{\mathbf{i}}^u$ which is not in $K_{\alpha}^u$ has its backwards iterates eventually in $K_{\alpha}^s$ and, hence, eventually expanded. Since the tangent vectors to $W_{\mathbf{i}}^u$ are eventually contracted in the past, they must be in $K_{\alpha}^u$.

We now return to the affine representatives $\tilde{f}_{i_j}$ of the maps $f_{i_j}$.

To obtain $g_{\mathbf{i}}$ satisfying (33), (34), it is convenient to use graph transform techniques as in [7], [12].

In view of Lemmas 4.3 and 6.3, we may assume that

$$K_0 > c, \ \epsilon_{ij}(w) < \frac{1}{4}, \ \epsilon_{22}(w) < \frac{1}{8} \tag{35}$$

for $w \in B_j$ where $c > 0$ is arbitrary. In the present section, it suffices to take $c > 4$. In section 8 below, we will take $c > 117$.

We define some function spaces.

Recall that $z_0 = \pi(\mathbf{i}), z_i = \pi(\sigma^i \mathbf{i})$ for $i \leq 0$. Let $\epsilon_i = \epsilon(\pi\sigma^i\mathbf{i}), B_i^u = B_{\epsilon_i}^u(z_i), B_i^s = B_{\epsilon_i}^s(z_i)$.

Let $\mathcal{G}_{0i}$ be the space of Lipschitz functions $g$ from $B_i^u$ to $B_i^s$ with Lipschitz constant less than or equal to 1. For such a $g$, let $graph(g) = \{(x, y) : y = g(x) \text{ for } x \in B_i^u\}$.

For $g_1, g_2 \in \mathcal{G}_{0i}$, set

$$d_{0i}(g_1, g_2) = \sup_{x \in B_i^u} \mid g_1 x - g_2 x \mid .$$

Let $\mathcal{G}_{1i}$ be the set of continuous functions $H : B_i^u \times \mathbf{R} \to \mathbf{R}$ such that for each $x \in B_i^u$, the map

$$v \to H(x, v)$$

is linear of norm no larger than 1.

Define the metric $d_{1i}$ on $\mathcal{G}_{1i}$ by

$$d_{1i}(H_1, H_2) = \sup_{x \in B_i^u, \mid v \mid \leq 1} \mid H_1(x, v) - H_2(x, v) \mid$$

Let $\mathcal{G}_{2i}$ be the set of continuous functions $J : B_i^u \times \mathbf{R} \times \mathbf{R} \to \mathbf{R}$ such that for each $x \in B_i^u$, the map

$$(v, w) \to J(x, v, w)$$

is symmetric and bilinear of norm no larger than $K_2$ for some constant $K_2$ to be specified later.

Set

$$d_{2i}(J_1, J_2) = \sup_{x \in B_i^u, \mid v \mid \leq 1} \mid J_1(x, v, v) - J_2(x, v, v) \mid$$

The spaces $(\mathcal{G}_{0i}, d_{0i}), (\mathcal{G}_{1i}, d_{1i}), (\mathcal{G}_{2i}, d_{2i})$ are bounded complete metric spaces.

Let $\mathbf{Z}^- = \{k \leq 0\}$ be the non-positive integers and consider the spaces

$$\mathcal{L}_0 = \{\phi : \mathbf{Z}^- \to \bigcup_i \mathcal{G}_{0i} : \phi_i \in \mathcal{G}_{0i} \ \forall \ i\}$$

$$\mathcal{L}_1 = \{\phi : \mathbf{Z}^- \to \bigcup_i \mathcal{G}_{1i} : \phi_i \in \mathcal{G}_{1i} \ \forall \ i\}$$

$$\mathcal{L}_2 = \{\phi : \mathbf{Z}^- \to \bigcup_i \mathcal{G}_{2i} : \phi_i \in \mathcal{G}_{2i} \ \forall \ i\}$$

with the metrics

$$\bar{d}_i(\phi, \psi) = \sum_{k\geq 0} \frac{1}{2^k} d_{ik}(\phi_k, \psi_k)$$

where $\phi, \psi \in \mathcal{L}_i, i = 0, 1, 2$.

The spaces $\mathcal{L}_i$ are also bounded complete metric spaces.

Let us recall the graph transform operator [7]. Let $f = \tilde{f}_{i_j}$ for some $i_j$, and let $g \in \mathcal{G}_{0j}$. Write $f(x, y) = (f_1(x, y), f_2(x, y))$, and let $(1, g) : B^u_j \to B^u_j \times B^s_j$ be the graph map defined by $(1, g)x = (x, gx)$.

We define

$$\Gamma(f, g) = f_2 \circ (1, g) \circ [f_1 \circ (1, g)]^{-1}.$$

It follows from our hyperbolicity assumptions the $\Gamma(f, g)$ is a well-defined mapping from $\mathcal{G}_{0j}$ to $\mathcal{G}_{0,j+1}$ for $j \leq -1$.

Returning now to the spaces $\mathcal{L}_i$ of sequences of functions, let us use the notation $\mathbf{g} = (g_k)_{k\leq 0}$, for elements of $\mathcal{L}_0, H = (H_k)_{k\leq 0}$, for elements of $\mathcal{L}_1$, and $J = (J_k)_{k\leq 0}$, for elements of $\mathcal{L}_2$. If $\mathbf{g} = (g_k)_{k\leq 0}$ is a sequence of $C^2$ functions, we write $D\mathbf{g} = (Dg_k)_{k\leq 0}, D^2\mathbf{g} = (D^2 g_k)_{k\leq 0}$.

We will define continuous maps $\Phi_0 : \mathcal{L}_0 \to \mathcal{L}_0, \Phi_1 : \mathcal{L}_0 \times \mathcal{L}_1 \to \mathcal{L}_1, \Phi_2 : \mathcal{L}_0 \times \mathcal{L}_1 \times \mathcal{L}_2 \to \mathcal{L}_2$, $\Xi_1 : \mathcal{L}_0 \times \mathcal{L}_1 \to \mathcal{L}_0 \times \mathcal{L}_1, \Xi_2 : \mathcal{L}_0 \times \mathcal{L}_1 \times \mathcal{L}_2 \to \mathcal{L}_0 \times \mathcal{L}_1 \times \mathcal{L}_2$ with the following properties.

FB1. $\Xi_1(g, H) = (\Phi_0(g), \Phi_1(g, H))$ and $\Xi_2(g, H, J) = (\Phi_0(g), \Phi_1(g, H), \Phi_2(g, H, J))$ for each $(g, H, J) \in \mathcal{L}_0 \times \mathcal{L}_1 \times \mathcal{L}_2$.

FB2. If $(g_k)_{k\leq 0}$ is a sequence of $C^2$ maps with $g_k \in \mathcal{G}_{0k}, Dg_k \in \mathcal{G}_{1k}, D^2 g_k \in \mathcal{G}_{2k}$ for all $k$, then $\Xi_2(\mathbf{g}, D\mathbf{g}, D^2\mathbf{g})_k = (\Gamma(\tilde{f}_{i_{k-1}}, g_{k-1}), D\Gamma(\tilde{f}_{i_{k-1}}, g_{k-1}), D^2\Gamma(\tilde{f}_{i_{k-1}}, g_{k-1}))$.

FB3. $\Phi_0$ is a contraction mapping; i.e., it is Lipschitz with Lipschitz constant less than 1.

FB4. The map $\Xi_1$ is a fiber contraction map over $\Phi_0$ in the sense of [7].

FB5. The map $\Xi_2$ is a fiber contraction map over $\Xi_1$.

Once these properties are established, we proceed as follows.

Let $z_0 = (x_0, y_0) \in \pi(\mathbf{i})$, let $\pi_2(x, y) = y$, and let $\mathbf{g} = (g_k)_{k\leq 0}$ be the sequence of constant maps

$$
\begin{aligned}
g_0(x) &= y_0 \\
g_{k-1}(x) &= \pi_2(\tilde{f}_{i_{-1}} \circ \ldots \circ \tilde{f}_{i_{k-1}})^{-1}(z_0).
\end{aligned}
$$

for $x \in B^u_{k-1}$.

Using the fiber contraction theorem of [7] we have that the sequence $\Xi_2^n(\mathbf{g}, D\mathbf{g}, D^2\mathbf{g}),\ n \geq 1$, converges to a fixed point $(\tilde{\mathbf{g}}, \tilde{H}, \tilde{J})$ of $\Xi_2$. Letting $\pi_0 : \mathcal{L}_0 \times \mathcal{L}_1 \times \mathcal{L}_2 \to \mathcal{L}_0, \pi_1 : \mathcal{L}_0 \times \mathcal{L}_1 \times \mathcal{L}_2 \to \mathcal{L}_1, \pi_2 : \mathcal{L}_0 \times \mathcal{L}_1 \times \mathcal{L}_2 \to \mathcal{L}_2$ be the natural projections, the definitions give

$$
\begin{aligned}
\pi_0 \Xi_2^n(\mathbf{g})_0 &= \Gamma(\tilde{f}_{i_{-1}} \circ \ldots \circ \tilde{f}_{i_{-n}}, g_{-n}) \\
\pi_1 \Xi_2^n(\mathbf{g})_0 &= D\Gamma(\tilde{f}_{i_{-1}} \circ \ldots \circ \tilde{f}_{i_{-n}}, g_{-n}) \\
\pi_2 \Xi_2^n(\mathbf{g})_0 &= D^2\Gamma(\tilde{f}_{i_{-1}} \circ \ldots \circ \tilde{f}_{i_{-n}}, g_{-n})
\end{aligned}
$$

Since all three of these sequences converge, it follows that

$$
\lim_{n\to\infty} \Gamma(\tilde{f}_{i_{-1}} \circ \ldots \circ \tilde{f}_{i_{-n}}, g_{-n}) = g_{\mathbf{i}}
$$

is $C^2$ with $Dg_{\mathbf{i}} = \lim \pi_1 \Xi_2^n(\mathbf{g})_0$ and $D^2 g_{\mathbf{i}} = \lim \pi_2 \Xi_2^n(\mathbf{g})_0$. The function $g_{\mathbf{i}}$ will be the $C^2$ function whose graph is contained in (and hence equals) $W^u_{\mathbf{i}}$.

Let us now define the maps $\Phi_i$ and establish their properties.

Let $f = \tilde{f}_{i_j}$ for some $i_j$ and let $g$ be a $C^2$ function such that $g \in \mathcal{G}_{0j}, Dg \in \mathcal{G}_{1j}, D^2 g \in \mathcal{G}_{2j}$.

Write $u(x) = [f_1 \circ (1, g)]^{-1}(x)$.

Then, differentiating $\Gamma(f, g) = f_2 \circ (1, g) \circ ([f_1 \circ (1, g)]^{-1}$ we get

$$
\begin{aligned}
D\Gamma(f, g) &= f_{2x}(ux, gux)Du(x) + f_{2y}(ux, gux)Dg(ux)Du(x) \\
&= f_{2x}Du(x) + f_{2y}Dg(ux)Du(x)
\end{aligned}
$$

$$
\begin{aligned}
D^2\Gamma(f, g) = \; & f_{2xx}Du(x)Du(x) + f_{2xy}Du(x)Dg(ux)Du(x) + f_{2yx}Du(x)Dg(ux)Du(x) \\
& + f_{2yy}Dg(ux)Dg(ux)Du(x)Du(x) + f_{2x}D^2u(x) \\
& + f_{2y}D^2g(ux)Du(x)Du(x) + f_{2y}Dg(ux)D^2u(x) \qquad (36)
\end{aligned}
$$

We can compute formulas for $Du, D^2u$ in terms of $f, g$ by differentiating the formula $f_1(ux, gux) = x$ twice and solving for $Du, D^2u$.

We get

$$Du(x) = \Big[f_{1x}(ux, gux) + f_{1y}(ux, gux)Dg(ux)\Big]^{-1} \tag{37}$$

and

$$\begin{aligned} D^2u(x) &= -Du(x)\Big[f_{1xx}(Du(x))^2 + 2f_{1xy}Dg(ux)Du(x)Du(x) \\ &\quad + f_{1yy}(Dg(ux))^2 Du(x)^2 + f_{1y}(Du)^2 D^2g(ux)\Big]. \end{aligned} \tag{38}$$

For $H \in \mathcal{G}_{1j}$, $J \in \mathcal{G}_{2j}$, let us write $H_x$ for the map $H(x, \cdot)$, $J_x$ for the map $J(x, \cdot, \cdot)$.

Define

$$D_1 = D_1(u, H)_x = \Big[f_{1x}(ux, gux) + f_{1y}(ux, gux)H_{ux}\Big]^{-1}$$

$$D_2(u, H, J)_x = -D_1\Big[f_{1xx}D_1D_1 + 2f_{1xy}H_{ux}D_1D_1 + f_{1yy}H_{ux}H_{ux}D_1D_1 + f_{1y}D_1D_1J_{ux}\Big]$$

$$R_1(f, g, h)_x = \Big[f_{2x}(ux, gux)) + f_{2y}(ux, gux)H_{ux}\Big]D_1$$

$$\begin{aligned} R_2(f, g, H, J)_x &= f_{2xx}D_1D_1 + f_{2xy}D_1H_{ux}D_1 + f_{2yx}D_1H_{ux}D_1 \\ &\quad + f_{2yy}H_{ux}H_{ux}D_1D_1 + f_{2x}D_2 + f_{2y}J_{ux}D_1D_1 + f_{2y}H_{ux}D_2 \end{aligned}$$

Finally, if $\mathbf{g} = (g_k)_{k\le 0} \in \mathcal{L}_0$, $H = (H_k)_{k\le 0} \in \mathcal{L}_1$, $J = (J_k)_{k\le 0} \in \mathcal{L}_2$, set

$$\Phi_0(\mathbf{g})_k = \Gamma(\tilde{f}_{i_{k-1}}, g_{k-1})$$

$$\Phi_1(\mathbf{g}, H)_k = (\Phi_0(\mathbf{g})_k, R_1(f_{i_{k-1}}, g_{k-1}, H_{k-1}))$$

$$\Phi_2(\mathbf{g}, H, J)_k = (\Phi_0(\mathbf{g})_k, R_1(\tilde{f}_{i_{k-1}}, g_{k-1}, H_{k-1}), R_2(\tilde{f}_{i_{k-1}}, g_{k-1}, H_{k-1}, J_{k-1})).$$

and define $\Xi_1, \Xi_2$ as in FB1.

Then, $\Xi_i$ and $\Phi_i$ satisfy properties FB1 and FB2 above.

Let us verify the fiber contraction properties of $\Xi_1, \Xi_2$.

**Fiber contraction property of $\Xi_1$.**

We first show that for fixed $f, g$ with $f = \tilde{f}_{i_j}$ and $g : B_j^u \to B_j^s$ a given Lipschitz map of Lipschitz constant no larger than 1, $R_1(f, g, \cdot)$ maps $\mathcal{G}_{1j}$ into $\mathcal{G}_{1,j+1}$ and is a contraction.

Since the graph of $g$ is in $Q$, the $C^0$ size of $\Gamma(f, g)$ is no larger than 1. This, and the overflowing property of $f$ on $B_j$ gives that $\Gamma(f, g)$ is a map from $B_{j+1}^u$ to $B_{j+1}^s$.

Let $Lip(\psi)$ be the Lipschitz constant of a map $\psi$.

As above, let $u(x) = [f_1 \circ (1, g)]^{-1}$.

Then,

$$Lip(u) \leq \frac{1}{\mid f_{1x} \mid (1 - \epsilon_{12})}$$

Using $\Gamma(f, g) = f_2 \circ (1, g) \circ [f_1 \circ (1, g)]^{-1}$, and the fact that $Lip(g) \leq 1$, we get

$$\begin{aligned} Lip(\Gamma(f, g)) &\leq Lip((f_2 \circ (1, g))Lip(u) \\ &\leq (\mid f_{2x} \mid + \mid f_{2y} \mid) Lip(u) \\ &\leq \frac{\epsilon_{21}}{1 - \epsilon_{12}} + \frac{\epsilon_{22}}{1 - \epsilon_{12}} \\ &\leq \frac{\epsilon_{21} + \epsilon_{22}}{1 - \epsilon_{12}} \leq 1 \end{aligned}$$

by (35). Thus, $\Gamma(f, g) \in \mathcal{G}_{1,j+1}$.

If $H, \tilde{H} \in \mathcal{G}_{1j}$, we have

$$\begin{aligned} \mid R_1(f, g, H) - R_1(f, g, \tilde{H}) \mid &\leq \mid (f_{2x} + f_{2y} H_{ux}) D_1(u, H) - (f_{2x} + f_{2y} \tilde{H}_{ux}) D_1(u, \tilde{H}) \mid \\ &\leq \mid f_{2y} \mid\mid D_1(u, H) \mid\mid H - \tilde{H} \mid \\ &\quad + (\mid f_{2x} \mid + \mid f_{2y} \mid\mid \tilde{H}_{ux} \mid) \mid D_1(u, H) - D_1(u, \tilde{H}) \mid \\ &\leq \frac{\epsilon_{22}}{1 - \epsilon_{12}} \mid H - \tilde{H} \mid \\ &\quad + (\mid f_{2x} \mid + \mid f_{2y} \mid\mid \tilde{H}_{ux} \mid) \mid D_1(u, H) - D_1(u, \tilde{H}) \mid \end{aligned}$$

To compute $| D_1(u,H) - D_1(u,\tilde{H}) |$, we use the formula

$$| G_1^{-1} - G_2^{-1} | \leq | G_1^{-1} || G_2^{-1} || G_1 - G_2 |$$

which follows immediately from the formula

$$| G_2^{-1} G_2 G_1^{-1} - G_2^{-1} G_1 G_1^{-1} | \leq | G_2^{-1} || G_2 - G_1 || G_1^{-1} |$$

Thus,

$$\begin{aligned} | D_1(u,H) - D_1(u,\tilde{H}) | &\leq | D_1(u,H) || D_1(u,\tilde{H}) || f_{1y} || H - \tilde{H} | \\ &\leq \frac{\epsilon_{12}}{| f_{1x} |(1-\epsilon_{12})^2} | H - \tilde{H} | \end{aligned}$$

Putting the above inequalities together, and using the fact that $| \tilde{H} | \leq 1$, we get

$$| R_1(f,g,H) - R_1(f,g,\tilde{H}) | \leq \frac{\epsilon_{22}}{1-\epsilon_{12}} | H - \tilde{H} | + \left[ \frac{\epsilon_{21}\epsilon_{12}}{(1-\epsilon_{12})^2} + \frac{\epsilon_{22}\epsilon_{12}}{(1-\epsilon_{12})^2} \right] | H - \tilde{H} |$$

Now, the fact that $R_1$ contracts the fibers follows from the estimates for $\epsilon_{ij}$ already given above.

The fiber norm of $R_2(f,g,H,J)$ and fiber contractions of $R_2(f,g,H,J)$ are obtained in the same way. We just write down the final estimates and leave the computations to the reader.

We have

$$\begin{aligned} | R_2(f,g,H,J) | &\leq \frac{4| D^2 f |}{| f_{1x} |^2 (1-\epsilon_{12})^2} + \frac{(\epsilon_{21}+\epsilon_{22}) 4 | D^2 f |}{(1-\epsilon_{12})^3 | f_{1x} |^2} \\ &\quad + \left[ \frac{(\epsilon_{21}+\epsilon_{22})\epsilon_{12}}{(1-\epsilon_{12})^3 | f_{1x} |} + \frac{\epsilon_{22}}{(1-\epsilon_{12})^2 | f_{1x} |} \right] | J | \\ &\leq \tilde{A}_1 + \tilde{A}_2 | J | \end{aligned}$$

and

$$| R_2(f,g,H,J) - R_2(f,g,H,\tilde{J}) | \leq \frac{(\epsilon_{21}+\epsilon_{22})\epsilon_{12}}{(1-\epsilon_{12})^3 | f_{1x} |} | J - \tilde{J} | + \frac{\epsilon_{22}}{(1-\epsilon_{12})^2 | f_{1x} |} | J - \tilde{J} |$$

Let us summarize the conditions we need to get the required properties of $R_1, R_2$.

$$\frac{\epsilon_{21}+\epsilon_{22}}{1-\epsilon_{12}} < 1 \tag{39}$$

$$\frac{\epsilon_{22}}{1-\epsilon_{12}} + \frac{\epsilon_{21}\epsilon_{12}}{(1-\epsilon_{12})^2} + \frac{\epsilon_{22}\epsilon_{12}}{(1-\epsilon_{12})^2} < 1 \tag{40}$$

$$\frac{(\epsilon_{21}+\epsilon_{22})\epsilon_{12}}{(1-\epsilon_{12})^3 | f_{1x} |} + \frac{\epsilon_{22}}{(1-\epsilon_{12})^2 | f_{1x} |} < 1 \tag{41}$$

Since $\epsilon_{12} < \frac{1}{4}$ and $K_0 > 4$, inequalities (39), (40), and (41) hold. Also, $| \tilde{A}_2 | < 1$. So, if we let $\tilde{K} > \frac{1}{1-\tilde{A}_2}$ and $K_2 = \tilde{K}\tilde{A}_1$, we have

$$| J | \leq K_2 \Rightarrow | R_2(f,g,H,J) | \leq K_2.$$

Hence, this $K_2$ is sufficient to define the space $\mathcal{G}_{2j}$.

**Proof of continuous dependence of the unstable manifolds $W^u_{\mathbf{i}}$ on the itineraries i.**

We have already noted that it suffices to prove that the functions $g_{\mathbf{i}}$ depend $C^2$ continuously on **i**.

It is clear that the maps $(f,g) \to \Gamma(f,g), (f,g,H) \to R_1(f,g,H), (f,g,H,J) \to R_2(f,g,H,J)$ are continuous. Since the spaces $\mathcal{G}_{0j}, \mathcal{G}_{1j}, \mathcal{G}_{2j}$ are bounded and the metrics on $\mathcal{L}_0, \mathcal{L}_1, \mathcal{L}_2$ give the product topologies, it follows that the maps $\Phi_0, \Xi_1, \Xi_2$ are continuous. Also our previous estimates give that, using the non-positive itineraries **i** as parameters, the family $(\Phi_0)_{\mathbf{i}}$ is a uniform family of contractions. Similarly, the families $(\Phi_1)_{\mathbf{i}}, (\Phi_2)_{\mathbf{i}}$ are families of uniform fiber contractions. Thus, the continuous dependence of $g_{\mathbf{i}}$ (and hence $W^u_{\mathbf{i}}$) follows from Propositions 5.1 and 5.2.

# 7 Fluctuation of Derivatives

We need to estimate quotients of the form

$$\frac{|\ D(f_{i_1} \circ \ldots \circ f_{i_n})_z(v_z)\ |}{|\ D(f_{i_1} \circ \ldots \circ f_{i_n})_w(v_w)\ |} \tag{42}$$

where $z, w$ are in a $K^u_\alpha$ curve $\gamma$ and $v_z, v_w$ are the unit tangent vectors to $\gamma$ at $z, w$, respectively.

The domains of the compositions $f_{i_1} \circ \ldots \circ f_{i_n}$ become narrow and possibly very non-convex. Since we wish to use the Mean Value Theorem in these domains, it will be convenient to choose certain star-shaped subdomains. This will be done in the next section. Here we present a useful Lemma.

Recall that a set $E$ is star-shaped relative to a point $z \in E$ if for any $w \in E$, the line segment joining $z$ to $w$ lies in $E$.

For a point $z \in E$ let $\delta_z(E)$ denote the diameter of the intersection of the horizontal line through $z$ and $E$.

Writing $f$ for one of the compositions above, assume that $Df$ maps the cone $K^u_\alpha$, into itself, expands it by at least $K_0 > 1$, and that $Df^{-1}$ maps the cone $K^s_\alpha$ into itself and expands it by at least $K_0$ as well.

For a subset $E$ of the domain of $f$ and $z \in E$, define

$$\Theta_z(f, E) = \sup_{w \in E} \frac{|\ D^2 f(w)\ |}{|\ f_{1x}(w)\ |} \delta_z(E)$$

where

$$|\ D^2 f(w)\ | = \max\{|\ f_{ijk}(w)\ | : i = 1, 2\ (j, k) = (x, x), (x, y), (y, y)\}.$$

**Lemma 7.1** *Let $E$ be a subset of the domain of $f$ which contains $z$ and is star-shaped relative to $z$. Let $\gamma$ be a $C^2$ curve in $E$ parametrized in the form $\gamma : x \to (x, g(x))$ where $g$ is a $C^2$ function such that $|\ Dg(x)\ | \leq \alpha$ and $|\ D^2 g(x)\ | \leq K_3$ for all $x$. Suppose $z, w \in \gamma, w \in E$, and $v_z, v_w$ are the unit tangent vectors to $\gamma$ at $z, w$, respectively. Let $\Theta = \Theta_z(f, E)$ and $\delta = \delta_z(E)$.*

*Then, there is a constant $C = C(\alpha, K_3) > 0$ such that*

$$\frac{|\ Df_z(v_z)\ |}{|\ Df_w(v_w)\ |} \leq exp\left(Cexp(C\Theta)\frac{|\ z - w\ |}{\delta}\right). \tag{43}$$

**Proof.**

We use the max norm $\mid (v_1, v_2) \mid = \max(\mid v_1 \mid, \mid v_2 \mid)$.

Let $z = (x, g(x)), w = (y, g(y)), v_z = (v_{1z}, v_{2z}), v_w = (v_{1w}, v_{2w})$.

Then, $\mid v_{1z} \mid = \mid v_{1w} \mid = 1$. Also, since $Df_z(v_z), Df_w(v_w)$ are in the cone $K^u_\alpha$, we have

$$\mid f_{1x}(z)v_{1z} + f_{1y}(z)v_{2z} \mid = \mid Df_z(v_z) \mid$$

and

$$\mid f_{1x}(w)v_{1w} + f_{1y}(w)v_{2w} \mid = \mid Df_w(v_w) \mid$$

So,

$$\begin{aligned}
\mid Df_w(v_w) \mid &= \mid f_{1x}(w)v_{1w} + f_{1y}(w)v_{2w} \mid \\
&= \mid f_{1x}(w))v_{1w} \mid \left(1 - \frac{\mid f_{1y}(w)v_{2w} \mid}{\mid f_{1x}(w)v_{1w} \mid}\right) \\
&\geq \mid f_{1x}(w) \mid (1 - \alpha^2)
\end{aligned}$$

and

$$\begin{aligned}
\frac{\mid Df_z(v_z) \mid}{\mid Df_w(v_w) \mid} &= 1 + \frac{\mid Df_z(v_z) \mid - \mid Df_w(v_w) \mid}{\mid Df_w(v_w) \mid} && (44) \\
&\leq exp(\frac{\mid Df_z(v_z) - Df_w(v_w) \mid}{\mid Df_w(v_w) \mid}) && (45) \\
&\leq exp(\frac{A_1(z, w) + A_2(z, w)}{(1 - \alpha^2)}) && (46)
\end{aligned}$$

where

$$A_1(z, w) = \frac{\mid Df_z(v_z) - Df_z(v_w) \mid}{\mid f_{1x}(w) \mid} \quad (47)$$

and

$$A_2(z, w) = \frac{\mid Df_z(v_w) - Df_w(v_w) \mid}{\mid f_{1x}(w) \mid} \quad (48)$$

We consider the two terms $A_1(z,w)$ and $A_2(z,w)$ separately.

We have

$$\begin{aligned} \mid Df_z(v_z) \; & - \; Df_z(v_w) \mid \\ & = \; \max(\mid f_{1x}v_{1z} + f_{1y}v_{2z} - f_{1x}v_{1w} - f_{1y}v_{2w} \mid, \mid f_{2x}v_{1z} + f_{2y}v_{2z} - f_{2x}v_{1w} - f_{2y}v_{2w} \mid) \end{aligned}$$

where the partial derivatives are all evaluated at $z$.

From Lemma 4.1 an upper bound for this last quantity is

$$\mid f_{1x}(z) \mid (1 + 2\alpha + \frac{1}{K_0^2} + \alpha^2) \mid v_z - v_w \mid$$

and this gives

$$A_1(z,w) \leq C(\alpha)\frac{\mid f_{1x}(z) \mid}{\mid f_{1x}(w) \mid} \mid v_z - v_w \mid$$

Now, $\mid v_z - v_w \mid$ is bounded above by the product of the maximum curvature of $\gamma$ and $\mid z - w \mid$. An upper bound for the curvature is the quantity $K_3$.

Let us use $C = C(\alpha, K_3)$ for possibly different values of $C$ below.

As in the proof of Lemma 4.2, we get

$$\frac{\mid f_{1x}(z) \mid}{\mid f_{1x}(w) \mid} \leq exp(C\Theta)$$

So,

$$A_1(z,w) \leq Cexp(C\Theta) \mid z - w \mid \leq C \exp(C\Theta)\frac{\mid z - w \mid}{\delta}. \tag{49}$$

Proceding similarly, the numerator of $A_2(z,w)$ is bounded above by

$$\begin{aligned} \max_{i=1,2}(\mid f_{ix}(z) - f_{ix}(w) \mid\mid v_{1w} \mid \; & + \; \mid f_{iy}(z) - f_{iy}(w) \mid\mid v_{2w} \mid) \\ & \leq \; 2 \max_{i=1,2,j=x,y} \mid f_{ij}(z) - f_{ij}(w) \mid \end{aligned}$$

Now,

$$| f_{ix}(z) - f_{ix}(w) | \leq | f_{ixx}(\tau) || z - w | + | f_{ixy}(\tau) || z - w |$$

and

$$| f_{iy}(z) - f_{iy}(w) | \leq | f_{iyx}(\tau_1) || z - w | + | f_{iyy}(\tau_1) || z - w |$$

for suitable $\tau, \tau_1$
which implies that

$$A_2(z,w) \leq C\Theta exp(C\Theta)\frac{| z - w |}{\delta}. \tag{50}$$

Using $C\Theta \leq exp(C\Theta)$, (49), (50) and a different $C$, we see that the proof of Lemma 7.1 is complete. QED

# 8 Distortion for compositions

In view of Lemma 7.1, to estimate quotients of the form (42), we will need to control the distortions of the compositions $\Theta_z(f_{i_1} \circ \ldots \circ f_{i_n})$ on appropriate sets.

Let $\mathbf{i} \in \Sigma$, and let $z \in W^s_{loc}(\pi\mathbf{i})$ be a point in the local stable manifold of $\pi(\mathbf{i})$. Write $i_j = i_j(z)$ for the $j-$th entry in the itinerary of $z$, and write $F^n(z) = f_{i_{n-1}} \circ \ldots f_{i_1} \circ f_{i_0}(z)$ so that $F^n(z) \in E_{i_n(z)}$ for all $n$.

For a curve $\gamma$, and $z, w \in \gamma$, let $v_z, v_w$ denote the unit tangent vectors to $\gamma$ at $z, w$, respectively.

As in section 2, let

$$E_{i_0 \ldots i_n} = E_{i_0} \bigcap f_{i_0}^{-1}(E_{i_1 \ldots i_n})$$

**Proposition 8.1** *(Bounded distortion of compositions) There is a constant $K_4 > 0$ such that for any $\mathbf{i} \in \Sigma$, any full width $K^u_\alpha$ curve $\gamma$ in $E_{i_0}$, and any $n > 0$, we have*

$$\frac{\mid DF^n_z(v_z) \mid}{\mid DF^n_w(v_w) \mid} \leq K_4 \tag{51}$$

*for any $z, w \in E_{i_0 \ldots i_n} \bigcap \gamma$.*

To prove this proposition, it will be convenient to cover the images $F^j(\gamma \bigcap E_{i_0 \ldots i_n})$ by small parallelograms in which the distortions $\Theta(F)$ become small, and to make use of affine coordinates as in section 6.

Let $E^s_z$ be the tangent space to $W^s_{loc}(\pi\mathbf{i})$ at $z$, and let $E^u_z$ be the tangent space to $\gamma$ at $z$. Writing $z_j$ for $F^j z,\ j \geq 0$, we translate these subspaces along the forward orbit of $z$ by defining

$$E^s_{z_j} = DF^j_z(E^s_z), E^u_{z_j} = DF^j(E^u_z),\ j \geq 0$$

This gives us a splitting of $T\mathbf{R}^2$ along the forward orbit of $z$ and the angles between the subspaces $E^s_{z_j}, E^u_{z_j}$ are uniformly bounded away from 0 by a constant that depends on $\alpha$.

Using these splittings, we can define affine coordinates along the forward orbit of $z$, giving local coordinate representatives $\tilde{f}_{i_j}$ of $f_{i_j}$, and small parallelograms $B_j = B^u_j \times B^s_j$ with sides parallel to the subspaces $E^u_{z_j}, E^s_{z_j}$ satisfying conditions analogous to those in B1, B2 following Lemma 6.3. As we have already noted, in view of Lemmas 4.3 and 6.3, we also can arrange for the conditions (35) to hold where $c > 117$.

In these affine coordinates, the subspaces $E^u_{z_j}, E^s_{z_j}$ become horizontal and vertical, repectively. As in section 6 we use the max norm in these coordinates, so each small $\epsilon-$ball $B_\epsilon(z_j) = B(z_j, \epsilon)$ will be a square of side length $2\epsilon$ centered at $z_j$.

If $E$ is any subset of $B_j$, and $z \in E$, let $C(z, E)$ denote the connected component of $E$ containing $z$. As in section 6, we may assume that

$$B_j \subset \bigcup_{w \in E_{i_j}} B(w, \bar{K}\delta_w(E_{i_j}))$$

where $\bar{K} > 0$ is a fixed constant.

For the remainder of this section we identify $f_{i_j}$ with its local coordinate representative $\tilde{f}_{i_j}$.

Thus, we may assume, for $w \in B_j$,

$$\mid f_{i_j 1x}(w) \mid \geq K_0 > 117 \tag{52}$$

$$\frac{\mid D^2 f_{i_j}(w) \mid}{\mid f_{i_j 1x}(w) \mid} \delta_{z_j}(E_{i_j}) < C_0 \tag{53}$$

$$\max(\epsilon_{12}(w), \epsilon_{21}(w)) < \epsilon_0, \ \epsilon_{22}(w) < \epsilon_0 \tag{54}$$

$$\bar{K}\delta_{z_j}(E_{i_j}) > diam(B_j) > C_1 \delta_{z_j}(E_{i_j}) \tag{55}$$

where $C_0, C_1, \epsilon_0$ are positive constants, $C_1 < C_0$, and $\epsilon_0 < \frac{1}{4}$.

We also may assume that $\gamma_j \equiv C(z_j, F^j\gamma \bigcap B_j)$ is a $K^u_{\epsilon_0}$ curve in $B_j$.

Let $\epsilon_1 \in (0, min(\frac{C_1}{2}, 1))$ be small enough so that

$$exp(156\epsilon_1 C_0)\frac{16}{15} < 2 \tag{56}$$

Let $B_{j,\epsilon_1} = B_j \bigcap B(z_j, \frac{\epsilon_1}{2}\delta_{z_j}(E_{i_j}))$.

The definition of $B_{j,\epsilon_1}$ implies that

$$\Theta_{z_j}(f_{i_j}, E) \leq \epsilon_1 C_0$$

for any subset $E \subset B_{j,\epsilon_1}$.

We use $\partial B$ to denote the boundary of a set $B$.

Since $f_{i_j}$ maps $E_{i_j}$ to a full-width rectangle in $Q$, there is a constant $K > 0$ such that

$$\delta_{z_j}(E_{i_j}) > K|\, f_{i_j 1x}(z_j)\,|^{-1}$$

Therefore, since $\delta_{z_j}(B_{j,\epsilon_1}) = \epsilon_1 \delta_{z_j}(E_{i_j})$, Lemma 4.2 provides a constant $K_5 > 0$ such that

$$dist(f_{i_j}(z_j), \partial f_{i_j}(\gamma_j \bigcap B_{j,\epsilon_1})) \geq K_5 \epsilon_1 \tag{57}$$

For $z_j \in E_{i_j}$, let

$$\tilde{B}_j = \begin{cases} B_{j,\epsilon_1} & \text{if } \frac{1}{2}\epsilon_1 \delta_{z_j}(E_{i_j}) < \frac{K_5 \epsilon_1}{2K_0} \\ B(z_j, \frac{K_5\epsilon_1}{2K_0}) & \text{if } \frac{1}{2}\epsilon_1 \delta_{z_j}(E_{i_j}) \geq \frac{K_5 \epsilon_1}{2K_0} \end{cases}$$

Thus, each $\tilde{B}_j \subseteq B_{j,\epsilon_1}$.

Since $f_{i_j}$ expands horizontal distances by at least $K_0$, we have that

$$dist(f_{i_j}(z_j), \partial f_{i_j}(\gamma_j \cap \tilde{B}_j)) \geq \frac{K_5 \epsilon_1}{2}$$

so

$$f_{i_j}(C(z_j, \gamma_j \cap \tilde{B}_j)) \supset C(f_{i_j} z_j, f_{i_j} \gamma_j \cap \tilde{B}_{j+1})$$

The set $\tilde{B}_{j,n} = \tilde{B}_{i_j} \bigcap F^{-1} \tilde{B}_{i_{j+1}} \bigcap \ldots \bigcap F^{-(n-1-j)} \tilde{B}_{i_{n-1}}$ is a narrow curvilinear rectangle around $z_j$.

Let

$$\alpha_{j,n} = dist(z_j, \gamma_j \bigcap \partial \tilde{B}_{j,n})$$

Let $E_{j,n}$ be the curvilinear rectangle whose left and right boundary curves are pieces of the left and right boundaries of $\tilde{B}_{j,n}$ and whose top and bottom boundary curves are horizontal line segments each of whose distance from $z_j$ is $\alpha_{j,n}$.

**Lemma 8.2** *The curvilinear rectangle $E_{j,n}$ is star-shaped relative to $z_j$.*

**Proof.** Let $\mathbf{b}_1, \mathbf{b}_2$ denote the left and right boundary curves of $E_{j,n}$ and let $\ell_1, \ell_2$ denote the top and bottom boundary curves (which are horizontal line segments).

Let $w \in E_{j,n}$, let $\ell_{j,w}$ denote the line segment joining $z_j$ to $w$, and let $\partial_{vert} E_{j,n}$ denote the union $\mathbf{b}_1 \bigcup \mathbf{b}_2$.

Since $z_j, w$ lie between the horizontal lines through $\ell_1, \ell_2$, any intersection of $\ell_{j,w}$ and the boundary of $E_{j,n}$, must be in $\partial_{vert} E_{j,n}$. Thus, to show that $\ell_{j,w}$ is contained in $E_{j,n}$ it suffices to show that

$$(\ell_{j,w} \setminus \{z_j, w\}) \cap \partial_{vert} E_{j,n} = \emptyset. \tag{58}$$

Assuming (58) fails we will get a contradiction.

By construction, $\mathbf{b}_1, \mathbf{b}_2$ are $K^s_{\epsilon_0}-$curves. This and the assumption that $\epsilon_0 < \frac{1}{4}$, imply that any line segment joining $z_j$ to a point $\bar{z}$ in $\partial_{vert} E_{j,n}$ must have slope no larger than $\frac{4}{3}$.

But since $z_j$ and $w$ lie in $E_{j,n}$, if (58) fails there is a point $\bar{z} \in \partial_{vert} E_{j,n}$ such that the line $\ell_{j,w}$ is parallel to the tangent vector to $\partial_{vert} E_{j,n}$ at $\bar{z}$. However, $\mathbf{b}_1, \mathbf{b}_2$ are $K^s_{\epsilon_0}-$curves, so their tangent vectors have slope no smaller than 4 which is our contradiction, proving Lemma 8.2.

**Lemma 8.3** *Fix $j \geq 0$. Then, for each $n > j$, we have*

$$\Theta_{z_j}(F^{n-j}, E_{j,n}) \leq 13 \epsilon_1 C_0. \tag{59}$$

**Proof.**

The proof is by induction on $n - j$. Clearly (59) holds for $n - j = 1$. Assuming it holds for $n - j$, we show it holds for $n + 1 - j$.

Let $z = z_j,\ f = F,\ g = F^{n-j},\ E_f = E_{n,n+1} = \tilde{B}_{gz},\ E_g = E_{j,n},\ h = f \circ g$, and $E_h = E_{j,n+1}$.

Let $\Delta f = \delta_{z_n}(E_f), \Delta g = \delta_{z_j}(E_g), \Delta h = \delta_{z_j}(E_h)$.

We use $\Theta(f) = \Theta_{gz}(f, E_f), \Theta(g) = \Theta_z(g, E_g), \Theta(h) = \Theta_z(h, E_h)$.

Consider the quotient

$$\frac{g_{1x}(w)}{g_{1x}(z)}$$

where $w \in E_h$.

Since the left and right boundary curves of $E_h$ are $K^s_{\epsilon_0}$ curves, and $\epsilon_0 < \frac{1}{4}$ we have that $|\ w - z\ | \leq 3\Delta h$.

Since both $g_{1x}(z), g_{1x}(w)$ have absolute value greater than 1, they have the same sign. Replacing $g$ by $-g$ if necessary, we may assume these signs are positive.

By the mean value theorem,

$$\begin{aligned}\mid \log g_{1x}(w) - \log g_{1x}(z) \mid &\leq \frac{\mid g_{1xx}(\tau) \mid}{\mid g_{1x}(\tau) \mid} \mid z - w \mid + \frac{\mid g_{1xy}(\tau) \mid}{\mid g_{1x}(\tau) \mid} \mid z - w \mid \\ &\leq 6\Theta(g)\frac{\Delta h}{\Delta g}\end{aligned}$$

so,

$$exp(-6\Theta(g)\frac{\Delta h}{\Delta g}) \leq \frac{\mid g_{1x}(w) \mid}{\mid g_{1x}(z) \mid} \leq exp(6\Theta(g)\frac{\Delta h}{\Delta g}). \tag{60}$$

Further, setting $\zeta = exp(12\Theta(g)\frac{\Delta_h}{\Delta_g})$, we have, for any $w, \tau \in E_h$,

$$\begin{aligned}\frac{\mid g_{1x}(w) \mid}{\mid g_{1x}(\tau) \mid} &= \frac{\mid g_{1x}(w) \mid}{\mid g_{1x}(z) \mid}\frac{\mid g_{1x}(z) \mid}{\mid g_{1x}(\tau) \mid} \\ &\leq \zeta. \end{aligned} \tag{61}$$

Similarly, if $\tau_1, \tau_2 \in E_g$, then

$$\frac{\mid g_{1x}(\tau_2) \mid}{\mid g_{1x}(\tau_1) \mid} \leq exp(12\Theta(g)) < exp(156\epsilon_1 C_0) < 2 \tag{62}$$

Also note that if $\ell_0$ is the full width horizontal line segment through $z$ in $E_h$, then $g(\ell_0)$ is a full width $K^u_{\epsilon_0}$ curve in $E_f$, and there is a $\tau \in \ell_0$ for which $\mid g_{1x}(\tau) \mid \Delta h = length(g(\ell_0)) \leq \frac{5}{4}\Delta f$.

This gives

$$\frac{\Delta h}{\Delta f} \leq \frac{5}{4\mid g_{1x}(\tau) \mid}. \tag{63}$$

Observe that it follows from the definition of $E_g$ and (57) that, for some $\tau_1 \in E_g$,

$$\mid g_{1x}(\tau_1) \mid \Delta g \geq K_5\epsilon_1.$$

So, by (62),

$$\begin{aligned}\frac{\Delta h}{\Delta g} &\leq \frac{5\Delta f}{4|\ g_{1x}(\tau)\ |}\frac{|\ g_{1x}(\tau_1)\ |}{K_5\epsilon_1} \\ &\leq \frac{5\Delta f}{2K_5\epsilon_1} \\ &\leq \frac{3}{K_0}\end{aligned} \tag{64}$$

Let us estimate

$$\Theta(h) = \max_{w\in E_h} \frac{|\ D^2 h(w)\ |}{|\ h_{1x}(w)\ |}\Delta h \tag{65}$$

Let

$$\eta = \frac{1}{1-\epsilon_{21}(g)\epsilon_{12}(f)} \leq \frac{16}{15}.$$

and

$$\zeta = exp(12\Theta(g)\frac{\Delta h}{\Delta g}) \leq exp(156\epsilon_1 C_0 \frac{3}{K_0})$$

Recall that $K_0$ and $\epsilon_1$ were chosen so that

$$K_0 > 117 \tag{66}$$

$$\eta\zeta < 2 \tag{67}$$

Recall the Chain Rule formulas (15)–(18).
By (15), we have

$$\begin{aligned}|\ h_{1x}(w)\ | &= |\ f_{1x}g_{1x} + f_{1y}g_{2x}\ | \\ &\geq |\ f_{1x}g_{1x}\ |(1-\epsilon_{21}(g)\epsilon_{12}(f)) \\ &= |\ f_{1x}g_{1x}\ |\eta^{-1}.\end{aligned}$$

Write $\epsilon_2(f) = max(\epsilon_{12}(f), \epsilon_{22}(f))$.

From (16) we get

$$\left|\frac{h_{ixx}(w)}{h_{1x}(w)}\Delta h\right| \leq \eta\left[\Theta(f)|\ g_{1x}(w)\ |\frac{\Delta h}{\Delta f}+2\Theta(f)|\ g_{2x}(w)\ |\frac{\Delta h}{\Delta f}+\Theta(f)\epsilon_{21}(g)|\ g_{2x}(w)\ |\frac{\Delta h}{\Delta f}\right.$$
$$\left.+\Theta(g)\max(1,\epsilon_{21}(f))(1+\epsilon_2(f))\frac{\Delta h}{\Delta g}\right]$$

Now, using (61), (63), (64), (66),(67),
we get

$$\begin{aligned}\left|\frac{h_{ixx}(w)}{h_{1x}(w)}\Delta h\right| &\leq \eta\left[\Theta(f)3\zeta(1+2\epsilon_{21}(g)+\epsilon_{21}(g)^2)+2\Theta(g)\frac{\Delta h}{\Delta g}\right]\\ &\leq 6\eta\zeta\Theta(f)+3\Theta(g)\frac{3}{K_0}\\ &\leq 12\Theta(f)+\frac{1}{13}\Theta(g)\end{aligned}$$

Similarly,

$$\begin{aligned}\left|\frac{h_{ixy}(w)}{h_{1x}(w)}\Delta h\right| &\leq \eta\,[3\Theta(f)\zeta(\epsilon_{12}(g)+\epsilon_{22}(g)+\\ &\quad\left.+\epsilon_{12}(g)\epsilon_{21}(g)+\epsilon_{22}(g)\epsilon_{21}(g))+2\Theta(g)\frac{\Delta h}{\Delta g}\right]\\ &\leq 3\Theta(f)+3\Theta(g)\frac{3}{K_0}\\ &\leq 3\Theta(f)+\frac{1}{13}\Theta(g)\end{aligned}$$

and

$$\begin{aligned}\left|\frac{h_{iyy}(w)}{h_{1x}(w)}\Delta h\right| &\leq \eta\left[3\Theta(f)\zeta(\epsilon_{12}(g)^2+2\epsilon_{22}(g)\epsilon_{12}(g)+\right.\\ &\quad\left.+\epsilon_{22}(g)^2\right]\\ &\leq 2\Theta(f)+3\Theta(g)\frac{3}{K_0}\\ &\leq 2\Theta(f)+\frac{1}{13}\Theta(g)\end{aligned}$$

In all cases we have

$$\Theta(h) \leq 12\Theta(f) + \frac{1}{13}\Theta(g) \leq 13\epsilon_1 C_0$$

proving Lemma 8.3.

**Proof of Proposition 8.1**.

The curvilinear rectangles $\tilde{B}_{j,n}$ are determined by the orbit segment $\{z_j = F^j(z)\}_{j=0}^n$. We write this as

$$\tilde{B}_{j,n} = \tilde{B}_{z,j,n}$$

There are analogous sets

$$\tilde{B}_{w,j,n} = \tilde{B}_{F^j w} \bigcap \ldots \bigcap F^{-(n-1-j)} \tilde{B}_{F^{n-1} w}$$

where $\tilde{B}_{F^\ell w}$ is a suitable small parallelogram centered at $F^\ell w$ for any $w \in E_{i_0 \ldots i_n} \cap \gamma$.

Let $\gamma, z, w, v_z, v_w$ be as in the hypotheses of the Proposition, and consider $F^n z, F^n w \in \gamma \cap E_{i_n}$.

We can connect these points by a chain $F^n z = F^n(w_1), F^n(w_2), \ldots F^n(w_k) = F^n(w)$ with $k \leq C_1(\alpha, \epsilon)$ such that, for every $\ell = 1, \ldots k-1$, and every $0 \leq j \leq n$,

$$F^j(w_{\ell+1}) \in \tilde{B}_{w_\ell, j, n}$$

then, it follows from Lemmas 7.1 and 8.3 that, for some constant $C_2(\alpha, \epsilon)$, we have

$$\frac{\mid DF^n_{w_\ell}(v_{w_\ell}) \mid}{\mid DF^n_{w_{\ell+1}}(v_{w_{\ell+1}}) \mid} \leq C_2(\alpha, \epsilon_1) \tag{68}$$

in the special affine coordinates centered at $w_\ell$. Changing back to the standard coordinates on $Q$ simply makes (68) hold with a different constant $C_2 = C_2(\alpha, \epsilon_1)$.

Then,

$$\frac{\mid DF^n_z(v_z) \mid}{\mid DF^n_w(v_w) \mid} = \prod_{\ell=1}^{k-1} \frac{\mid DF^n_{w_\ell}(v_{w_\ell}) \mid}{\mid DF^n_{w_{\ell+1}}(v_{w_{\ell+1}}) \mid} \leq C_2^{k-1}$$

proving Proposition 8.1.

## 9 Sinai Local Measures

For two points $z_1, z_2$ in an unstable manifold $W^u_{\mathbf{i}}$ and unit tangent vectors $v_1, v_2$ to $W^u_{\mathbf{i}}$ at $z_1, z_2$, respectively, let $D^u F(z_i) = \mid DF_{z_i}(v_i) \mid$ denote the Jacobian of $F$ at $z_i$ along $W^u_{\mathbf{i}}$. We know that $W^u_{\mathbf{i}}$ is a full-width $K^u_\alpha$ curve in $E_{i_0}$. Also, the curve $f_{i_0} W^u_{\mathbf{i}}$ is a full width $K^u_\alpha$ curve in $Q$.

**Proposition 9.1** *Suppose* $\mathbf{i} = (\dots i_{-n} \dots i_0)$ *is an arbitrary infinite non-positive itinerary and let* $W^u_{\mathbf{i}}$ *denote its unstable manifold. Write* $f_{i_0} W^u_{\mathbf{i}} =$ *graph* $g_{\mathbf{i}}$ *where* $g_{\mathbf{i}} : I \to I$ *is the* $C^2$ *function given in Theorem 6.1. Suppose* $x_1, x_2 \in I$ *and* $z_1 = (x_1, g_{\mathbf{i}} x_1), z_2 = (x_2, g_{\mathbf{i}} x_2)$*. Then, the infinite product*

$$\xi(x_1, x_2, \mathbf{i}) = \prod_{s=1}^{\infty} \frac{D^u F(F^{-s} z_1)}{D^u F(F^{-s} z_2)} \tag{69}$$

*converges and depends continuously on* $(x_1, x_2, \mathbf{i})$.

*Moreover, there is a constant* $K_6 > 0$ *independent of* $(x_1, x_2, \mathbf{i})$ *such that*

$$K_6^{-1} < \xi(x_1, x_2, \mathbf{i}) < K_6 \tag{70}$$

**Proof**.

It clearly suffices to prove the upper bound in (70) since interchanging $z_1$ and $z_2$ would then give the lower bound.

Let $\bar{z}_1 = f^{-1}_{i_0} z_1, \bar{z}_2 = f^{-1}_{i_0} z_2$ so that $\bar{z}_1, \bar{z}_2 \in W^u_{\mathbf{i}}$.

We use the local coordinates $\tilde{f}_{i_j}$ and rectangles $\tilde{B}_j$ of the previous section. To avoid confusion, we will use $\tilde{F}^s(z) = (\tilde{f}_{i_{-1}} \circ \dots \circ \tilde{f}_{i_{-s}})(z)$ and $\tilde{F}^{-s} = (\tilde{f}_{i_{-1}} \circ \dots \circ \tilde{f}_{i_{-s}})^{-1}$ instead of identifying $F, f_{i_j}$ with $\tilde{F}, \tilde{f}_{i_j}$ as in the preceding section. We use $\tilde{B}_z$ for the affine neighborhood centered at $z$.

In our local coordinates, with $z \in E_{i_0}, W^u_{\mathbf{i}} \bigcap \tilde{B}_z$ becomes a $K^u_{\epsilon_0}$ curve. Also, there is a sequence $\bar{z}_1 = w_1, \dots, w_k = \bar{z}_2$ of points in $W^u_{\mathbf{i}}$ such that

$$d(w_{j+1}, w_j) < \frac{\epsilon_1 K_5}{2} \tag{71}$$

and

$$k \leq \left[ \frac{2}{\epsilon_1 K_5} \right] + 1. \tag{72}$$

Recall that $\tilde{F}^{-1}w_j = \tilde{f}_{i_{-1}}^{-1}w_j$.

Further, $\tilde{f}_{i_{-1}}(\tilde{f}_{i_{-1}}^{-1}W^u_{\mathbf{i}} \cap \tilde{B}_{\tilde{F}^{-1}w_j})$ contains the intersection of $W^u_{\mathbf{i}}$ with the ball of radius $\frac{\epsilon_1 K_5}{2}$ about $w_j$. Since $w_{j+1}$ is in this latter set, we have $\tilde{F}^{-1}w_{j+1} \in \tilde{B}_{\tilde{F}^{-1}w_j}$.

Analogously, we have $\tilde{F}^{-s}w_{j+1} \in \tilde{B}_{\tilde{F}^{-s}w_j}$ for every $s \geq 1$.

Now, there is a constant $C_6 > 0$ such that

$$\begin{aligned}\prod_{s=1}^{\infty} \frac{D^uF(F^{-s}z_1)}{D^uF(F^{-s}z_2)} &\leq C_6 \prod_{s=1}^{\infty} \frac{D^uF(F^{-s}\bar{z}_1)}{D^uF(F^{-s}\bar{z}_2)} \\ &= \prod_{s=1}^{\infty}\prod_{j=1}^{k-1} \frac{D^uF(F^{-s}w_j)}{D^uF(F^{-s}w_{j+1})}\end{aligned}$$

so, to prove Proposition 9.1 it suffices to show

$$\prod_{s=1}^{\infty} \frac{D^uF(F^{-s}w_j)}{D^uF(F^{-s}w_{j+1})} \leq K_7 \tag{73}$$

for some $K_7 > 0$ and any $j$.

Since the angles between $E^u_z$ and $E^s_z$ are bounded by a constant depending on $\alpha$, the linear maps $D\tilde{F}(\tilde{F}^{-s}(w_j))$ and $DF(F^{-s}(w_j))$ are conjugate by a linear map whose images on unit vectors are bounded above and below by constants which depend only on $\alpha$. A similar statement holds replacing $w_j$ by $w_{j+1}$. Hence, there is a constant $C_5 = C_5(\alpha)$ such that, for any $s \geq 1$ and any $j$,

$$C_5^{-1}\frac{D^u\tilde{F}(\tilde{F}^{-s}w_j)}{D^u\tilde{F}(\tilde{F}^{-s}w_{j+1})} \leq \frac{D^uF(F^{-s}w_j)}{D^uF(F^{-s}w_{j+1})} \leq C_5\frac{D^u\tilde{F}(\tilde{F}^{-s}w_j)}{D^u\tilde{F}(\tilde{F}^{-s}w_{j+1})} \tag{74}$$

so, it suffices to find $K_7 > 0$ such that

$$\prod_{s=1}^{\infty} \frac{D^u\tilde{F}(\tilde{F}^{-s}w_j)}{D^u\tilde{F}(\tilde{F}^{-s}w_{j+1})} \leq K_7 \tag{75}$$

By Lemmas 7.1 and 8.3, there is a constant $K_1 > 0$ such that

$$\frac{\mid D\tilde{F}^N_z(v_z) \mid}{\mid D\tilde{F}^N_{\bar{z}}(v_{\bar{z}}) \mid} \leq K_1 \tag{76}$$

for any $N > 1, \bar{z} \in E_{z,N}$ and unit vectors $v_z, v_{\bar{z}}$ tangent to $W^u(z), W^u(\bar{z})$, respectively.

Let $N$ be large enough so that

$$\tau = \frac{K_1}{K_0^N K_5 \epsilon_1} < 1. \tag{77}$$

By definition, $\tilde{F}^{N-1}(E_{\tilde{F}^{-sN}w_j,N})$ is a full-width subrectangle of $\tilde{B}_{\tilde{F}^{N-1-sN}w_j}$. So, the $\tilde{F}$ image of a full-width horizontal line segment in $\tilde{B}_{\tilde{F}^{N-1-sN}w_j}$ contains a curve of horizontal width at least $\epsilon_1 K_5$.

Thus, setting $\tilde{w}_j = \tilde{F}^{-1}w_j$ and $\delta_{i_s} = \delta_{\tilde{F}^{-sN}\tilde{w}_j}(E_{\tilde{F}^{-sN}\tilde{w}_j,N})$, we have

$$\frac{\epsilon_1 K_5}{\mid D^u\tilde{F}^N(\tau_N)\mid} \le \delta_{i_s} \le \frac{1}{K_0^N} \tag{78}$$

for some $\tau_N \in W^u(\tilde{F}^{-sN}\tilde{w}_j) \cap E_{\tilde{F}^{-sN}\tilde{w}_j,N}$.

Then,

$$\begin{aligned}
\mid \tilde{F}^{-sN}\tilde{w}_{j+1} - \tilde{F}^{-sN}\tilde{w}_j \mid &= \mid \tilde{F}^{-sN+N}\tilde{w}_{j+1} - \tilde{F}^{-sN+N}\tilde{w}_j \mid \frac{1}{\mid D^u\tilde{F}^N(\tilde{\tau_N})\mid} \\
&\le \frac{\epsilon_1 K_5}{\mid D^u\tilde{F}^N(\tau_N)\mid}\frac{\mid D^u\tilde{F}^N(\tau_N)\mid}{\epsilon_1 K_5 \mid D^u\tilde{F}^N(\tilde{\tau_N})\mid}\mid \tilde{F}^{-sN+N}\tilde{w}_{j+1} - \tilde{F}^{-sN+N}\tilde{w}_j \mid \\
&\le \delta_{i_s}\frac{K_1}{\epsilon_1 K_5}\mid \tilde{F}^{-sN+N}\tilde{w}_{j+1} - \tilde{F}^{-sN+N}\tilde{w}_j \mid \\
&\le \delta_{i_s}\delta_{i_{s-1}}\dots\delta_{i_1}\left(\frac{K_1}{\epsilon_1 K_5}\right)^s \mid \tilde{w}_{j+1} - \tilde{w}_j \mid
\end{aligned}$$

giving

$$\begin{aligned}
\frac{\mid \tilde{F}^{-sN}\tilde{w}_{j+1} - \tilde{F}^{-sN}\tilde{w}_j \mid}{\delta_{i_s}} &\le \left(\frac{K_1}{\epsilon_1 K_5}\right)^s \frac{1}{K_0^{N(s-1)}}\mid \tilde{w}_{j+1} - \tilde{w}_j \mid \qquad (79) \\
&\le \frac{K_1}{\epsilon_1 K_5}\tau^{s-1}.
\end{aligned}$$

Hence,

$$\prod_{s=1}^{\infty} \frac{D^u\tilde{F}(\tilde{F}^{-s}w_j)}{D^u\tilde{F}(\tilde{F}^{-s}w_{j+1})} = \prod_{s=1}^{\infty} \frac{D^u\tilde{F}^N(\tilde{F}^{-sN-1}w_j)}{D^u\tilde{F}^N(\tilde{F}^{-sN-1}w_{j+1})}$$
$$= \prod_{s=1}^{\infty} \frac{D^u\tilde{F}^N(\tilde{F}^{-sN}\tilde{w}_j)}{D^u\tilde{F}^N(\tilde{F}^{-sN}\tilde{w}_{j+1})}$$
$$\leq exp(\sum_{s=1}^{\infty} C_7\tau^s) \equiv K_7$$

using Lemma 7.1 and (79).

Since the functions $g_{\mathbf{i}}$ depend continuously on $\mathbf{i}$ in the $C^2$ topology, the continuity statement in 9.1 follows from the fact that given $\epsilon > 0$, there is an $N_0 > 0$ such that if $N \geq N_0$, we have

$$\left| \prod_{j=N_0}^{\infty} \frac{D^uF(F^{-s}z_1)}{D^uF(F^{-s}z_2)} - 1 \right| < \epsilon$$

which is immediate from the proof just given.

This proves Proposition 9.1.

For a $C^2$ curve $\gamma$ in $Q$, let $\rho_\gamma$ denote the Riemannian measure on $\gamma$.

From Proposition 9.1 we get the existence of the following limit

$$\lim_{n\to\infty} \prod_{s=1}^{n} \frac{D^uF(F^{-s}z_1)}{D^uF(F^{-s}z_2)} = \xi(z_1, z_2) = \xi_{\mathbf{i}}(z_1, z_2) \tag{80}$$

for any two points $z_1, z_2 \in W^u_{\mathbf{i}}$. Letting $\gamma$ denote $W^u_{\mathbf{i}}$, we can use $\rho_\gamma$ and the ratios $\xi(z_1, z_2)$ obtained in the preceding limits to get special measures on the unstable manifolds. More precisely, following Sinai in [13], Lecture 16, we define

$$\nu_{z_1,\gamma}(A) = \int_A \xi(z_1, z_2) d\rho_\gamma(z_2).$$

It is easy to see that if $z_3$ is another point in $\gamma$, then $\nu_{z_3,\gamma}(A) = \xi(z_3, z_1)\nu_{z_1,\gamma}(A)$, so the measures $\nu_{z_1,\gamma}$ and $\nu_{z_3,\gamma}$ are simply rescalings of each other. In particular, if $A, B \subset \gamma$ and $\nu_{z_1,\gamma}(B) < \infty$, then $\frac{\nu_{z_1,\gamma}(A)}{\nu_{z_1,\gamma}(B)}$ is independent of $z_1$.

For $z_1 \in \gamma \bigcap \tilde{Q}$, let $E_{i_1}$ be the element of $\{E_i\}$ containing $Fz_1$, and let $\gamma_1 = W^u_{\sigma\mathbf{i}}$.

The family of measures $\{\nu_{z_1,\gamma}\}$ is invariant in the sense that if $A, B \subset \gamma$, $F(A), F(B) \subset \gamma_1$, and $\nu_{z_1,\gamma}(B) < \infty$, then $\frac{\nu_{z_1,\gamma}(A)}{\nu_{z_1,\gamma}(B)} = \frac{\nu_{Fz_1,\gamma_1}(FA)}{\nu_{Fz_1,\gamma_1}(FB)}$.

We call the family of measures $\nu_{z_1,\gamma}$ *Sinai local measures* or just *local measures.*

# 10 Absolute Continuity of the Stable Foliation

We know that for each non-negative itinerary $\mathbf{a} = (a_0, a_1, \ldots)$ there is a $C^1 K^s_\alpha$ curve $W^s(\mathbf{a}) = \bigcap_{n \geq 0} E_{a_0 \ldots a_{n-1}}$ of full height in $Q$.

Note that two points in $\tilde{Q}$ with different forward itineraries have disjoint stable manifolds since the interiors of the $E_i's$ are disjoint. Thus, the set $\{W^s(\mathbf{a}) : \pi\mathbf{a} \in \tilde{Q}\}$ is a foliation of its union. We call this the *stable foliation.* Let $\mathcal{W} = \{W^s(\mathbf{a}) : \pi\mathbf{a} \in \tilde{Q}\}$ denote this foliation. We denote the union $\bigcup\{W^s(\mathbf{a}) : \pi\mathbf{a} \in \tilde{Q}\}$ by $\mathcal{W}^+$. Note that $\mathcal{W}^+$ is a Borel subset of $Q$ of full two-dimensional Lebesgue measure in $Q$. For any two full width $C^2$ $K^u_\alpha$ curves $\gamma, \eta$, let $\pi_{\gamma\eta}$ be the holonomy projection from $\gamma$ to $\eta$ along the foliation $\mathcal{W}$. That is, for $z \in \gamma \bigcap \mathcal{W}^+$, and $W^s(\mathbf{a})$ the leaf of $\mathcal{W}$ which contains $z$, $\pi_{\gamma\eta}(z)$ is the unique point of intersection of $W^s(\mathbf{a})$ and $\eta$. As above, for any $C^2$ $K^u_\alpha-$curve $\gamma$, let $\rho_\gamma$ denote the Riemannian measure on $\gamma$. Recall that the foliation $\mathcal{W}$ is called *absolutely continuous* if

(AC-1) each full-width $C^2$ $K^u_\alpha$ curve $\gamma$ meets $\mathcal{W}^+$ in a set of positive $\rho_\gamma$ measure

and

(AC-2) the image measure $\pi_{\gamma\eta\star}\rho_\gamma$ is equivalent to the measure $\rho_\eta$.

**Proposition 10.1** *The foliation $\mathcal{W}$ is absolutely continuous*

Before we can prove Proposition 10.1, we need a couple of Lemmas.

The next Lemma is well-known and elementary. Since the proof is short, we include it for completeness.

**Lemma 10.2** *Suppose that $x_1, x_2, \ldots$ and $y_1, y_2, \ldots$ are sequences of numbers in the open unit interval $(0, 1)$ such that*

$$-\sum_{i \geq 1} x_i \log y_i < \infty. \tag{81}$$

*For $\epsilon > 0$ and non-ngative integer $n$, let $D_n = \{i : y_i \leq exp(-\epsilon n)\}$.*

*Then,*

$$\sum_{n\geq 1}\sum_{i\in D_n} x_i < \infty. \tag{82}$$

**Proof.**
For each $n \geq 1$, let

$$E_n = \{i : exp(-\epsilon(n+1)) < y_i \leq exp(-\epsilon n)\}$$

Then, $D_n = \bigsqcup_{j\geq n} E_j$, where $\bigsqcup$ denotes disjoint union, so

$$\sum_{n\geq 1}\sum_{i\in D_n} x_i = \sum_{n\geq 1}\sum_{j\geq n}\sum_{i\in E_j} x_i.$$

Letting $c_j = \sum_{i\in E_j} x_i$, this last sum is just

$$\begin{array}{cccccc} c_1 & + & c_2 & & + & c_3 + \ldots \\ & + & c_2 & & + & c_3 + \ldots \\ & & & & + & c_3 + \ldots \\ & & & & \vdots & \\ & = & \sum_{j\geq 1} jc_j & & & \end{array}$$

Now, $i \in E_j$ implies that $-\log y_i \geq \epsilon j$ or $-x_i \log y_i \geq \epsilon j x_i$ which gives

$$\begin{aligned} -\sum_{i\in E_j} x_i \log y_i &\geq \sum_{i\in E_j} \epsilon j x_i \\ &= \epsilon j c_j \end{aligned}$$

Hence,

$$\epsilon \sum_{j\geq 1} jc_j \leq \sum_{j\geq 1}\sum_{i\in E_j} -x_i \log y_i \leq -\sum x_i \log y_i < \infty$$

which implies that $\sum_{j\geq 1} jc_j < \infty$. QED.

In the next lemma, we will use the geometric condition $G3$.
Each $z \in \tilde{Q}$ has a unique forward itinerary $(a_0(z), a_1(z), \ldots)$ with $F^n(z) \in$ *int* $E_{a_n(z)}$.

**Lemma 10.3** *Let $\gamma$ be a $C^2$ $K^u_\alpha-$curve of full-width in $Q$ such that $\rho_\gamma(\gamma \bigcap \tilde{Q}) > 0$. Let $\epsilon > 0$. For $\rho_\gamma-$ almost all points $z \in \gamma \bigcap \tilde{Q}$, there is a positive integer $n(z) > 0$ such that if $n \geq n(z)$, then*

$$\delta_{F^n(z)}(E_{a_n(z)}) > exp(-\epsilon n).$$

**Proof.**

For ease of notation, if $A$ is a subset of $\gamma$, let us write $\mid A \mid$ for $\rho_\gamma(A)$.

Let $D_n = \{i \geq 1 : \delta_{i,min} < e^{-\epsilon n}\}$.

In view of lemma 10.2, the condition G3 implies that

$$\sum_{n \geq 1} \sum_{i \in D_n} \delta_{i,max} < \infty. \tag{83}$$

Let $V_n = \{z \in \gamma \bigcap \tilde{Q} : \delta_{F^n(z)}(E_{a_n(z)}) \leq e^{-\epsilon n}\}$.

We will show

$$\sum_{n \geq 1} \mid V_n \mid < \infty. \tag{84}$$

Once this is done, the Borel-Cantelli Lemma gives that $\rho_\gamma-$almost all points of $\gamma$ lie in at most finitely many of the $V_n's$ which proves Lemma 10.3.

Let $\mathcal{A}_n$ be the set of finite itineraries $(a_0, \ldots, a_{n-1})$ which occur for points in $\tilde{Q}$.

For a given finite sequence $a_0, a_1, \ldots, a_{n-1} \in \mathcal{A}_n$, let $V_n(a_0, \ldots, a_{n-1}) = \{z \in V_n : F^i z \in E_{a_i}$ for $0 \leq i < n\}$. Then,

$$V_n(a_0, \ldots, a_{n-1}) = \bigcup_{i \geq 0} \left(\gamma \bigcap E_{a_0 \ldots a_{n-1} i} \bigcap \tilde{Q}\right)$$

and this last union is disjoint.

Also, $V_n$ is the disjoint union of the $V_n(a_0, \ldots, a_{n-1})$ as these finite itineraries vary in $\mathcal{A}_n$.

The bounded distortion of compositions (Proposition 8.1) gives us a constant $K > 0$ such that for $(a_0, \ldots, a_{n-1} i) \in \mathcal{A}_{n+1}$, and $z \in \gamma \bigcap E_{a_0 \ldots a_{n-1} i}$,

$$\frac{\mid \gamma \bigcap E_{a_0 \ldots a_{n-1} i} \mid}{\mid \gamma \bigcap E_{a_0 \ldots a_{n-1}} \mid} \leq K \delta_{F^n(z)}(E_i)$$

Also, the definition of $V_n(a_0,\dots,a_{n-1})$ gives us that $\delta_{a_n(z),min} \le e^{-\epsilon n}$;i.e., that $a_n(z) \in D_n$.

Thus,

$$V_n(a_0,\dots,a_{n-1}) \subseteq \bigsqcup_{i\in D_n} \gamma \bigcap E_{a_0,\dots,a_{n-1}i} \bigcap \tilde{Q}$$

This gives

$$\begin{aligned}
\mid V_n \mid &\le \sum_{(a_0\dots a_{n-1})\in\mathcal{A}_n} \sum_{i\in D_n} \mid \gamma \cap E_{a_0,\dots,a_{n-1}i} \mid \\
&= \sum_{a_0\dots a_{n-1}} \sum_{i\in D_n} \frac{\mid \gamma \cap E_{a_0,\dots,a_{n-1}i} \mid}{\mid \gamma \cap E_{a_0\dots a_{n-1}} \mid} \mid \gamma \cap E_{a_0\dots a_{n-1}} \mid \\
&\le \sum_{a_0\dots a_{n-1}} \sum_{i\in D_n} K\delta_{i,max} \mid \gamma \cap E_{a_0\dots a_{n-1}} \mid \\
&\le \sum_{i\in D_n} K\delta_{i,max}
\end{aligned}$$

Hence, (84) is a consequence of (83). QED.

**Lemma 10.4** *For any full-width $K^u_\alpha$ curve $\gamma$,*

$$\rho_\gamma(\gamma \bigcap \tilde{Q}) = 1$$

**Proof.** The curve $\gamma$ cannot meet both the upper and lower boundaries of $Q$. For definiteness, we suppose that $\gamma$ does not meet the lower boundary of $Q$. The other case is similar.

Then, there are constants $\alpha_1 > 0, \alpha_2 > 0$ and a $C^1$ diffeomorphism $\phi$ from $Q$ onto a curvilinear subrectangle $Q_1$ of $Q$ such that

1. $\phi$ maps the upper boundary of $Q$ onto $\gamma$ and maps the lower boundary of $Q$ onto itself.
2. $D\phi(K^u_\alpha) \subset K^u_{\alpha_1}$
3. $D\phi^{-1}(K^s_\alpha) \subset K^s_{\alpha_2}$.

Let $\tilde{\gamma} = \phi^{-1}(\gamma)$ denote the upper boundary of $Q$.

Since a subset $A$ of $\gamma$ has full $\rho_{\gamma}$ measure if and only if $\phi^{-1}(A)$ has full $\rho_{\tilde{\gamma}}$ measure, it suffices to prove that

$$\rho_{\tilde{\gamma}}(\phi^{-1}(\gamma \cap \tilde{Q})) = 1 \tag{85}$$

Let $\tilde{E}_i = \phi^{-1}(E_i)$,

$$\tilde{\delta}_{i,max} = \max_{z \in \tilde{E}_i} \delta_z(\tilde{E}_i)$$

and

$$\tilde{\delta}_{i,min} = \min_{z \in \tilde{E}_i} \delta_z(\tilde{E}_i)$$

The properties of $\phi$ guarantee that

$$-\sum_i \tilde{\delta}_{i,max} \log \tilde{\delta}_{i,min} < \infty \tag{86}$$

Now, $Q_1 \cap \tilde{Q}$ has full Lebesgue measure in $Q_1$, so, $\phi^{-1}(\tilde{Q})$ has full Lebesgue measure in $Q$. Thus, for almost all horizontal lines $\ell$ in $Q$, we have that $\ell \cap \phi^{-1}(\tilde{Q})$ has full Riemannian measure.

To complete the proof of Lemma 10.4, we will prove that $\rho_\ell(\ell \bigcap \phi^{-1}(\tilde{Q}))$ varies continuously with $\ell$.

This is a consequence of the following.

For any $\epsilon > 0$, there is an $N = N(\epsilon) > 0$ such that for any horizontal full-width line segment $\ell$,

$$\rho_\ell(\ell \bigcap \bigcup_{i \geq N} \tilde{E}_i) < \epsilon$$

which is, in turn, a consequence of

$$\sum_{i \geq N} diam(\ell \bigcap \tilde{E}_i) < \epsilon.$$

Since the vertical boundaries of the $\tilde{E}_i's$ are $K^s_{\alpha_2}$ curves, there is a constant $C(\alpha_2) > 0$ such that, for all $i, m(\tilde{E}_i) > C(\alpha_2)(\tilde{\delta}_{i,max})^2$. So, $\tilde{\delta}_{i,max} \to 0$ as $i \to \infty$.

By (86) and Lemma 10.2 with $x_i = \tilde{\delta}_{i,max}, y_i = \tilde{\delta}_{i,min}$, given $\epsilon > 0$, we can find $n_0 > 0$ such that

$$\sum_{\tilde{\delta}_{i,min} < 2^{-n_0}} \tilde{\delta}_{i,max} < \epsilon$$

Now, take $N$ such that $i \geq N$ implies that $\tilde{\delta}_{i,min} < 2^{-n_0}$.
This gives $\sum_{i \geq N} diam(\ell \bigcap \tilde{E}_i) \leq \sum_{\tilde{\delta}_{i,min} < 2^{-n_0}} \tilde{\delta}_{i,max} < \epsilon$ as required . QED

For future use let us observe that the argument in the last proof actually works for all $K^u_\alpha$ curves uniformly to prove

**Lemma 10.5** *Given $\epsilon > 0$, there is an integer $N(\epsilon) > 0$ such that for every $K^u_\alpha$ curve $\gamma$, we have*

$$\rho_\gamma(\gamma \bigcap (\bigcup_{i \geq N} E_i)) < \epsilon$$

**Proof of Proposition 10.1.**
We use $\nu << \mu$ for $\nu$ is absolutely continuous with respect to $\mu$, and $\nu \sim \mu$ for $\nu << \mu$ and $\mu << \nu$.
Let $\gamma, \eta$ be two $C^2$ full-width $K^u_\alpha$ curves.
In what follows we restrict our measures to $\tilde{Q}$. Thus, when we write $\rho_\gamma(A)$ we mean $\rho_\gamma(A \bigcap \tilde{Q})$.
We will show that

$$\pi_{\gamma\eta\star}\rho_\gamma << \rho_\eta \tag{87}$$

Once this is done, interchanging $\gamma$ and $\eta$, we have $\pi_{\eta\gamma\star}\rho_\eta << \rho_\gamma$.
So, $\rho_\eta = \pi_{\gamma\eta\star}(\pi_{\eta\gamma\star}\rho_\eta) << \pi_{\gamma\eta\star}\rho_\gamma$ or $\rho_\eta \sim \pi_{\gamma\eta\star}\rho_\gamma$ as required for the proof of Proposition 10.1.
We know that $\rho_\gamma(\gamma \bigcap \tilde{Q}) = 1$.
Let $B \subset \gamma \bigcap \tilde{Q}$ be such that $\rho_\gamma(B) > 0$.
Let $K_1 \in (1, K_0)$.
By Lemma 10.3, for almost all $z \in \gamma$, there is an $n(z) > 0$ such that $n \geq n(z)$ implies

$$\delta_{F^n(z)}(E_{a_n(z)}) > K_1^{-n}. \tag{88}$$

From standard measure theory, we can take a compact set $A \subset B$ such that $\rho_\gamma(A) > 0$ and there is an $n(A) > 0$ such that (88) holds for all $n \geq n(A)$ and all $z \in A$.

We will show that there is a constant $K > 0$ such that

$$\rho_\eta(\pi_{\gamma_\eta}(A)) \geq K^{-1}\rho_\gamma(A) \tag{89}$$

This, in turn gives $\rho_\eta(\pi_{\gamma_\eta}(B)) > 0$ to prove (87).

Since $K_1 < K_0$, and $dist(F^n(z), F^n(\pi_{\gamma_\eta}(z))) \leq const \cdot K_0^{-n}$, we may assume that, for $z \in A$ and large $n$,

$$\frac{1}{2} < \frac{diam(F^n(E_{a_0(z)\ldots a_n(z)} \bigcap \gamma))}{diam(F^n(E_{a_0(z)\ldots a_n(z)} \bigcap \eta))} < 2 \tag{90}$$

For a unit vector $v$ tangent to the curve $\gamma$ at $\tau$ and a positive integer $n$, let us write $D_\gamma F^n(\tau)$ for $DF^n(v)$.

Now, for $z \in A$, there are points $\tau_n \in \gamma, \tilde{\tau} \in \eta$ such that

$$diam(\gamma \bigcap E_{a_0(z)\ldots a_n(z)}) \mid D_\gamma F^n(\tau_n) \mid = diam(F^n(\gamma \bigcap E_{a_0(z)\ldots a_n(z)}))$$

and

$$diam(\eta \bigcap E_{a_0(z)\ldots a_n(z)}) \mid D_\eta F^n(\tilde{\tau}_n) \mid = diam(F^n(\eta \bigcap E_{a_0(z)\ldots a_n(z)})).$$

We claim

(AC-3) there is a constant $K = K(A) > 0$ such that for all $z \in A$ and $n \geq 0$

$$K^{-1} < \frac{\mid D_\gamma F^n(z) \mid}{\mid D_\eta F^n(\pi_{\gamma_\eta}(z)) \mid} < K. \tag{91}$$

Assuming (AC-3) for the moment, we see that there is a possibly different $K > 0$ such that, for all $n \geq 0, z \in A$, we have

$$K^{-1} < \frac{diam(\gamma \cap E_{a_0(z)\ldots a_n(z)})}{diam(\eta \cap E_{a_0(z)\ldots a_n(z)})} < K. \tag{92}$$

But, for large $n$, as $z$ varies in $A$, the sets $\gamma \cap E_{a_0(z)\ldots a_n(z)}$ form a covering of $A$ by small intervals and the sets $\eta \cap E_{a_0(z)\ldots a_n(z)}$ form a covering of $\pi_{\gamma_\eta}(A)$ by small intervals. This gives (89) and concludes the proof of Proposition 10.1.

**Proof of (AC-3):**

Let $z_n = F^n(z), w_n = F^n(\pi_{\gamma\eta}(z))$ for each $n \geq 0$. We use affine coordinates centered as $z_n$ as in our earlier sections. We use the splitting $T_{z_n}\mathbf{R}^2 = E^u_{z_n} \oplus E^s_{z_n}$ in which $E^u_{z_n}$ contains $DF^n(v_z)$ and $E^s_{z_n}$ is tangent to $W^s_{loc}(z_n)$ at $z_n$.

Let $\tilde{F}$ denote the representative of $F$ in these coordinates, and let $\tilde{B}_n$ be the small parallelogram centered at $z_n$ as before. We may and do assume that $K_0 > 3$.

Write $v_{z_n}, v_{w_n}$ for the unit vectors tangent to $F^n(\gamma)$ at $z_n$ and $F^n(\eta)$ at $w_n$, respectively.

Now,

$$\frac{D_\gamma F^n(z)}{D_\eta F^n(\pi_{\gamma\eta}(z))} \leq const \cdot \prod_{s=1}^{n-1} \frac{\mid D\tilde{F}_{z_s}(v_{z_s}) \mid}{\mid D\tilde{F}_{w_s}(v_{w_s}) \mid} \tag{93}$$

so, it suffices to show

$$\frac{\mid D\tilde{F}_{z_n}(v_{z_n}) \mid}{\mid D\tilde{F}_{w_n}(v_{w_n}) \mid} \leq exp(a_n) \tag{94}$$

where

$$\sum_{n\geq 1} a_n < const \cdot \log K \tag{95}$$

to prove (AC-3).

Write $\delta_n$ for $\delta_{F^n(z)}(E_{a_n(z)})$.

In our affine coordinates, $v_{z_n} = \begin{pmatrix} 1 \\ 0 \end{pmatrix}$.

Since $dist(F^n(\pi_{\gamma\eta}(z), F^n(z))$ is exponentially smaller than $\delta_n$ for large $n$ and $\mid v_{w_n} - v_{z_n} \mid \to 0$ as $n \to \infty$, there is an $n_0 = n_0(A)$ such that $n \geq n_0$ implies $w_n \in \tilde{B}_n$ and $v_{w_n} \in K^u_{\epsilon_0}$. (Here $\epsilon_0 < \frac{1}{4}$ as in section 6).

Below, we use various constants $C_s$, $1 \leq s \leq 8$, which are independent of $n$ and $z \in A$ and are defined in the first equation in which they appear.

As in the proof of lemma 7.1,

$$\frac{\mid D\tilde{F}_{z_n}(v_{z_n}) \mid}{\mid D\tilde{F}_{w_n}(v_{w_n}) \mid} \leq exp(A_{1,n} + A_{2,n})$$

where

$$A_{1,n} \leq C_1 | \, v_{z_n} - v_{w_n} \, |$$

and

$$A_{2,n} \leq C_2 \frac{| \, z_n - w_n \, |}{\delta_n}.$$

Since

$$\frac{| \, z_n - w_n \, |}{\delta_n} \leq C_3 \left( \frac{K_1}{K_0} \right)^n$$

it suffices to show

$$| \, v_{w_n} - v_{z_n} \, | \leq C_4 \left( \frac{K_1}{K_0} \right)^{n-1} \tag{96}$$

for all $n$ to prove (94), (95), and (AC-3).
Writing $D\tilde{F}_{w_{n-1}}(v_{w_{n-1}}) = (\xi_n, \eta_n)$ and $v_{w_n} = (u_n^1, u_n^2)$ we have

$$\begin{aligned} \xi_n &= \tilde{F}_{1x}(w_{n-1}) u_{n-1}^1 + \tilde{F}_{1y}(w_{n-1}) u_{n-1}^2 \\ \eta_n &= \tilde{F}_{2x}(w_{n-1}) u_{n-1}^1 + \tilde{F}_{2y}(w_{n-1}) u_{n-1}^2 \end{aligned}$$

and

$$| \, D\tilde{F}_{w_{n-1}}(v_{w_{n-1}}) \, | = | \, \xi_n \, |.$$

Thus,

$$u_n^1 = 1, \; u_n^2 = \frac{\eta_n}{| \, \xi_n \, |}.$$

This gives

$$\begin{aligned} | \, v_{w_n} - v_{z_n} \, | &\leq \frac{| \, \eta_n \, |}{| \, \xi_n \, |} \\ &\leq \left( \frac{1}{1 - \epsilon_0^2} \right) \left( \frac{| \, \tilde{F}_{2x}(w_{n-1}) \, |}{| \, \tilde{F}_{1x}(w_{n-1}) \, |} + \frac{| \, \tilde{F}_{2y}(w_{n-1}) \, |}{| \, \tilde{F}_{1x}(w_{n-1}) \, |} | \, u_{n-1}^2 \, | \right) \end{aligned}$$

Using $\tilde{F}_{2x}(z_{n-1}) = 0$, we get

$$\begin{aligned}
\frac{| \tilde{F}_{2x}(w_{n-1}) |}{| \tilde{F}_{1x}(w_{n-1}) |} &\leq \frac{| \tilde{F}_{2xx}(\tau) |}{| \tilde{F}_{1x}(w_{n-1}) |} | w_{n-1} - z_{n-1} | + \frac{| \tilde{F}_{2xy}(\tau) |}{| \tilde{F}_{1x}(w_{n-1}) |} | w_{n-1} - z_{n-1} | \\
&\leq C_5 \frac{| w_{n-1} - z_{n-1} |}{\delta_{n-1}}
\end{aligned}$$

for suitable $\tau$.
Analogously,

$$\frac{| \tilde{F}_{2y}(w_{n-1}) |}{| \tilde{F}_{1x}(w_{n-1}) |} \leq \frac{| \tilde{F}_{2y}(z_{n-1}) |}{| \tilde{F}_{1x}(w_{n-1}) |} + C_6 \frac{| w_{n-1} - z_{n-1} |}{\delta_{n-1}}$$

which gives

$$\begin{aligned}
\frac{| \eta_n |}{| \xi_n |} &\leq \frac{1}{(1-\epsilon_0^2)K_0^2} | u_{n-1}^2 | + C_7 \frac{| w_{n-1} - z_{n-1} |}{\delta_{n-1}} \\
&\leq \frac{1}{(1-\epsilon_0^2)K_0^2} | u_{n-1}^2 | + C_8 \left(\frac{K_1}{K_0}\right)^{n-1}
\end{aligned}$$

Inductively, we assume

$$| u_{n-1}^2 | \leq 2C_8 \left(\frac{K_1}{K_0}\right)^{n-2}$$

and get

$$| u_n^2 | = \frac{| \eta_n |}{| \xi_n |} \leq \frac{2C_8}{(1-\epsilon_0^2)K_0^2} \left(\frac{K_1}{K_0}\right)^{n-2} + C_8 \left(\frac{K_1}{K_0}\right)^{n-1}$$

Since $\frac{3}{K_0} < 1 < K_1$ and $\epsilon_0 < \frac{1}{4}$, we get $\frac{2}{(1-\epsilon_0^2)K_0^2} < \frac{K_1}{K_0}$ and

$$| u_n^2 | \leq 2C_8 \left(\frac{K_1}{K_0}\right)^{n-1}$$

which proves (96).

# 11 Construction of an SRB measure

We wish to use a construction analogous to that of Sinai in [13] to construct our SRB measure. There are several difficulties which appear.

1. The family of unstable manifolds $\{W^u(z)\}$ does not form a measurable partition of the attractor $\Lambda$ in $Q$.
2. The underlying set $\Lambda$ is not compact, so care has to exercised in the taking of limits of iterates of measures.

We will see that these problems can be handled by lifting the required construction to the symbolic space $\Sigma$, getting a measure there, compactifying, getting a limit measure which is supported on $\Sigma$, and projecting back into $Q$.

We have defined a continuous map $\pi$ from $\Sigma$ into $Q$ as follows. For $\mathbf{a} \in \Sigma$ with $\mathbf{a} = (\ldots a_{-1}a_0a_1\ldots)$,

$$\{\pi(\mathbf{a})\} = \bigcap_{n\geq 0} E_{a_0\ldots a_n} \cap f_{a_0}^{-1} S_{a_{-n}\ldots a_0}$$

Let $\sigma$ be the left shift automorphism on $\Sigma$. For each $\mathbf{a} \in \Sigma$, we have local stable and unstable sets defined by

$$W^s_{loc}(\mathbf{a}) = \{\mathbf{b} : a_i = b_i, i \geq 0\}$$

$$W^u_{loc}(\mathbf{a}) = \{\mathbf{b} : a_i = b_i, i \leq 0\}$$

We have the local stable and unstable sets in $Q$ as well:

$$W^u_{loc}(\pi\mathbf{a}) = \bigcap_{n\geq 0} f_{a_0}^{-1} S_{a_{-n}\ldots a_0}$$

$$W^s_{loc}(\pi\mathbf{a}) = \bigcap_{n\geq 0} E_{a_0\ldots a_n}$$

Each $W^u_{loc}(\pi\mathbf{a})$ is a $K^u_\alpha$ curve which has full width in $E_{a_0}$, and each $W^s_{loc}(\pi\mathbf{a})$ is a $K^s_\alpha$ curve of full height in $Q$.

Note that if $\mathbf{a} = (\ldots a_i \ldots), \mathbf{b} = (\ldots b_i \ldots), \pi\mathbf{a}, \pi\mathbf{b} \in \tilde{Q}$, and $a_i \neq b_i$ for some $i \geq 0$, then $W^s_{loc}(\pi\mathbf{a}) \bigcap W^s_{loc}(\pi\mathbf{b}) = \emptyset$.

Thus, the map $\pi : W^u_{loc}(\mathbf{a}) \bigcap \pi^{-1}\tilde{Q} \to W^u_{loc}(\pi\mathbf{a}) \bigcap \tilde{Q}$ is a one-to-one, continuous onto map for each $\mathbf{a} \in \pi^{-1}(\tilde{Q})$. By standard results, it is a Borel isomorphism.

Recall the functions $\xi(z_1, z_2)$ and the Sinai local measures $\nu_{z_1,\gamma}$ defined at the end of section 9.

We now use them to define finite measures on the local unstable sets $W^u_{loc}(\mathbf{a})$ in $\Sigma$.

Write $\tilde{W}^u_{loc}(\mathbf{a}) = W^u_{loc}(\mathbf{a}) \bigcap \pi^{-1}\tilde{Q}$.

If $\gamma = W^u_{loc}(\pi\mathbf{a})$, then $\gamma \bigcap \tilde{Q}$ has full Riemannian measure in $\gamma$, and the Borel isomorphism $\pi : \tilde{W}^u_{loc}(\mathbf{a}) \to W^u_{loc}(\pi\mathbf{a}) \bigcap \tilde{Q}$ allows us to transfer the Riemannian measure $\rho_\gamma$ from $\gamma \bigcap \tilde{Q}$ up to $\tilde{W}^u_{loc}(\mathbf{a})$. We call this measure $\rho_{\mathbf{a}}$. It clearly only depends on the non-positive indices of $\mathbf{a}$.

For $z, w \in \tilde{W}^u_{loc}(\mathbf{a})$, let

$$\bar{\xi}(z, w) = \xi(\pi z, \pi w)$$

where $\xi(\cdot, \cdot)$ is the density of the Sinai local measure defined at the end of Section 9.

Next, for $z \in \tilde{W}^u_{loc}(\mathbf{a})$, we define a finite measure $\nu_z$ on $\tilde{W}^u_{loc}(\mathbf{a})$ by

$$\nu_z(A) = \int_A \bar{\xi}(z, w) d\rho_{\mathbf{a}}(w)$$

These measures have the following properties

1. For $z_1, z_2 \in \tilde{W}^u_{loc}(\mathbf{a})$, and $A \subset \tilde{W}^u_{loc}(\mathbf{a})$

$$\nu_{z_1}(A) = \bar{\xi}(z_1, z_2)\nu_{z_2}(A)$$

2. If $A, B \subset \tilde{W}^u_{loc}(\mathbf{a}), z_1 \in \tilde{W}^u_{loc}(\mathbf{a}), \nu_{z_1}(B) > 0$, and $\sigma(A), \sigma(B) \subset \tilde{W}^u_{loc}(\sigma\mathbf{a})$, then $\nu_{\sigma z_1}(\sigma B) > 0$ and

$$\frac{\nu_{\sigma z_1}(\sigma A)}{\nu_{\sigma z_1}(\sigma B)} = \frac{\nu_{z_1}(A)}{\nu_{z_1}(B)}$$

It follows from these facts that if $\nu_{z_1}(B) > 0$, for some $z_1$, then $\nu_{z_2}(B) > 0$ for any $z_2$, and the normalized measure $\nu_B(A) = \frac{\nu_{z_1}(A \bigcap B)}{\nu_{z_1}(B)}$ is independent of the choice of $z_1 \in \tilde{W}^u_{loc}(\mathbf{a})$. Moreover, the normalized measures are $\sigma-$invariant in the following sense: if $A$ and $B$ are as in 2 above, then $\sigma_\star \nu_B(\sigma(A)) = \nu_{\sigma(B)}(\sigma(A))$. We will call the measures $\nu_z$, *local measures* or *Sinai measures.*

For a point $\mathbf{a} \in \Sigma$, with local unstable set $\tilde{W}^u_{loc}(\mathbf{a})$, let $\nu_{\mathbf{a},norm}$ be its normalized local measure. Thus,

$$\nu_{\mathbf{a},norm}(A) = \frac{\nu(A)}{\nu_{\mathbf{a}}(\tilde{W}^u_{loc}(\mathbf{a}))}$$

for every $A \subset \tilde{W}^u_{loc}(\mathbf{a})$.

For each $i \geq 1$, let $V_i = \{\mathbf{a} \in \Sigma : a_0 = i\}$, and fix a local stable set $S_i \subset V_i$. Thus, $S_i = W^s_{loc}(z_i)$ where $z_i$ is a particular point in $V_i$. Let $\mathcal{M}_i$ be the partition of $V_i$ into local unstable sets. The quotient set $V_i/\mathcal{M}_i$ is in one-to-one correspondence with $S_i$, so the partition $\mathcal{M}_i$ is measurable with respect to any complete Borel probability measure on $V_i$. Let $\mathcal{M} = \bigcup_i \mathcal{M}_i$. Since $\Sigma$ is a countable disjoint union of the $V_i's$, $\mathcal{M}$ is a measurable partition of $\Sigma$ for any complete Borel probability measure.

For convenience, we will say that a Borel partition $\mathcal{M}$ is *measurable* with respect to a Borel Probability measure $\mu$, if it is equal mod zero to a measurable Borel partition of the Borel completion of the measure $\mu$. This allows us to discuss systems of conditional measures, etc, with respect to arbitrary measurable Borel partitions of Borel probability measures.

Now fix an element $z_0 \in \pi^{-1}(\tilde{Q})$, and let $\tilde{W}^u_{loc}(z_0)$ be its local unstable set. Let $\nu_0$ be the associated normalized Sinai measure.

**Theorem 11.1** *The sequence of averages*

$$\nu_n = \frac{1}{n} \sum_{k=0}^{n-1} \sigma^k_\star \nu_0$$

*converges weakly to a measure $\bar{\mu}$ on $\Sigma$ which is $\sigma-$invariant, ergodic, and the conditional measures of $\bar{\mu}$ with respect to the partition $\mathcal{M}$ coincide with the normalized Sinai measures on elements of $\mathcal{M}$.*

The proof will require several steps.

Let $\mathbf{N}$ be the set of positive integers, and let $\bar{\mathbf{N}} = \mathbf{N} \bigcup \{\infty\}$ be its one-point compactification. We put a metric on $\bar{\mathbf{N}}$ making it isometric to $\{0, 1, \frac{1}{2}, \frac{1}{3} \ldots\} \subset \mathbf{R}$ with the standard metric. Let $\bar{\Sigma} = \bar{\mathbf{N}}^{\mathbf{Z}}$ with the product topology and let $\bar{\sigma} : \bar{\Sigma} \to \bar{\Sigma}$ be the shift. The set $\Sigma$ is a dense $\bar{\sigma}-$invariant subset of $\bar{\Sigma}$.

We take a subsequence $\{\nu_{n_k}\}$ of $\{\nu_n\}$ which converges to a measure $\bar{\mu}$ on $\bar{\Sigma}$.

**Claim 1**: The measure $\bar{\mu}$ is supported on $\Sigma$. That is,

$$\bar{\mu}(\bar{\Sigma} \setminus \Sigma) = 0.$$

**Proof.**

A point $\mathbf{a} \in \bar{\Sigma} \setminus \Sigma$ has $a_i = \infty$ for some $i$. Fixing $i$, let $J_i = \{\mathbf{a} \in \bar{\Sigma} : a_i = \infty\}$.

We will show that, given $\epsilon > 0$, there is an open neighborhood $U_i$ of $J_i \setminus \Sigma$ such that for all $n \geq 1 - i$,

$$\sigma_\star^n(\nu_0)(U_i) \leq \epsilon \tag{97}$$

This will imply that $\bar{\mu}(J_i \setminus \Sigma) = 0$. Since this holds for every $i$, Claim 1 follows.

Let $\epsilon_1 > 0$ be a small number to be chosen later.

From Lemma 10.5, there is an $N > 0$, such that for every $K^u_\alpha$ curve $\gamma$,

$$\rho_\gamma \left( \bigcup_{j \geq N} E_j \right) < \epsilon_1 \tag{98}$$

Let $\rho_0 = \rho_{z_0}$ be the lift to $\tilde{W}^u_{loc}(z_0)$ of the Riemannian measure on $W^u_{loc}(\pi z_0) \cap \tilde{Q}$, and let $\ell_0 = \rho_0(\tilde{W}^u_{loc}(z_0))$.

By Proposition 9.1, for any $E \subset \tilde{W}^u_{loc}(z_0)$,

$$K_6^{-2} \, \frac{\rho_0(E)}{\ell_0} \leq \nu_0(E) \leq K_6^2 \, \frac{\rho_0(E)}{\ell_0} \tag{99}$$

Given a non-negative itinerary $\mathbf{a} = (a_0 a_1 \ldots)$, let

$$V_{a_0 \ldots a_n} = \{\mathbf{b} \in \tilde{W}^u_{loc}(z_0) : b_i = a_i, i = 0, \ldots, n\}$$

By Proposition 8.1, for any $n \geq 1$, if $\gamma_n = F^n(W^u_{loc}(\pi z_0))$, then

$$K_4^{-1}\ \rho_{\gamma_n}(E_{a_n}) \leq \frac{\rho_0(V_{a_0\ldots a_n})}{\rho_0(V_{a_0\ldots a_{n-1}})} \leq K_4\ \rho_{\gamma_n}(E_{a_n}) \tag{100}$$

Setting $U_i = \{\mathbf{a} \in \bar{\Sigma} : a_i \geq N\}$, we see that (98) and (100) imply that, if $n + i \geq 1$, then

$$\rho_0(V_{a_0\ldots a_{n+i-1}} \cap \sigma^{-n}U_i) \leq K_4\epsilon_1\rho_0(V_{a_0\ldots a_{n+i-1}}) \tag{101}$$

Also, $\tilde{W}^u_{loc}(z_0) \cap \sigma^{-n}(U_i)$ is the disjoint union

$$\bigsqcup_{a_0\ldots a_{n+i-1}} V_{a_0\ldots a_{n+i-1}} \cap \sigma^{-n}U_i$$

So,

$$\begin{aligned}
(\sigma^n_\star\nu_0)(U_i) &= \nu_0(\sigma^{-n}(U_i)) \\
&= \nu_0(\tilde{W}^u_{loc}(z_0) \cap \sigma^{-n}U_i) \\
&= \sum_{a_0\ldots a_{n+i-1}} \nu_0(V_{a_0\ldots a_{n+i-1}} \cap \sigma^{-n}U_i) \\
&\leq \sum_{a_0\ldots a_{n+i-1}} \frac{K_6^2}{\ell_0}\rho_0(V_{a_0\ldots a_{n+i-1}} \cap \sigma^{-n}U_i) \\
&= \sum_{a_0\ldots a_{n+i-1}} \frac{K_6^2}{\ell_0}\frac{\rho_0(V_{a_0\ldots a_{n+i-1}} \cap \sigma^{-n}U_i)}{\rho_0(V_{a_0\ldots a_{n+i-1}})} \cdot \rho_0(V_{a_0\ldots a_{n+i-1}}) \\
&\leq \sum_{a_0\ldots a_{n+i-1}} \frac{K_6^2}{\ell_0}K_4\epsilon_1\rho_0(V_{a_0\ldots a_{n+i-1}}) \\
&= \epsilon_1\frac{K_6^2}{\ell_0}K_4\ell_0 \\
&= \epsilon_1 K_6^2 K_4
\end{aligned}$$

Hence, if we set $\epsilon_1 = \frac{\epsilon}{K_6^2K_4}$, we get (97), and Claim 1 is proved.

The measure $\bar{\mu}$ is clearly invariant under the shift $\sigma$.

We extend the partition $\mathcal{M}$ of $\Sigma$ to $\bar{\Sigma}$ by adding the element $\bar{\Sigma} \setminus \Sigma$. We will also use the letter $\mathcal{M}$ to denote this extended partition. We let $\bar{V}_i$ denote the closure of $V_i$ in $\bar{\Sigma}$, and let $\mathcal{M}_i$ denote the restriction of $\mathcal{M}$ to $\bar{V}_i$.

Let $\tilde{\pi} : \bar{\Sigma} \to \bar{\Sigma}/\mathcal{M}$ be the natural projection. Let $\tilde{\mu} = \tilde{\pi}_\star \bar{\mu}$ be the induced measure on $\bar{\Sigma}/\mathcal{M}$.

There is a system of conditional measures $\bar{\mu}_C$ on $C \in \mathcal{M}$ defined for $\tilde{\mu}$-almost all $C \in \mathcal{M}$.

**Claim 2:** For $\tilde{\mu}-$almost all $C$, $\bar{\mu}_C = \nu_C$.

**Proof.**

Let us use $\bar{A}$ for the closure of a subset $A \subset \bar{\Sigma}$ in $\bar{\Sigma}$.

Let $\phi : \bar{\Sigma} \to \mathbf{R}$ be a continuous function supported in $\bar{V}_i$ for some $i$.

For each $n \geq 0$, the measure $\sigma^n_\star \nu_0$ is supported on countably many $C's$ in $\mathcal{M}$, and these $C's$ are local unstable sets.

The conditional measure $(\sigma^n_\star \nu_0)_C$ is then just the restriction of $\sigma^n_\star \nu_0$ to $C$ normalized.

But, the invariance property of quotients of the Sinai measures gives, for $A \subset C$,

$$\begin{aligned} \frac{(\sigma^n_\star \nu_0)(A)}{(\sigma^n_\star \nu_0)(C)} &= \frac{\nu_0(\sigma^{-n}A)}{\nu_0(\sigma^{-n}C)} \\ &= \frac{\nu_C(A)}{\nu_C(C)} = \nu_C(A) \end{aligned}$$

Thus,

(*) the conditional measure $(\sigma^n_\star \nu_0)_C$ is equal to the normalized Sinai measure $\nu_C$ when $C$ is a local unstable set in $\sigma^n(\tilde{W}^u_{loc}(z_0))$.

and this implies

(**) the conditional measure $(\nu_n)_C$ equals $\nu_C$ on each local unstable set $C$ in $\Sigma$ such that $\nu_n(C) > 0$.

Let $S_i$ be the stable set of $z_i \in V_i$. Its closure $\bar{S}_i$ is the local stable set of $z_i$ in $\bar{\Sigma}$. This is a compact subset of $\bar{\Sigma}$ and may be identified with $\bar{V}_i/\mathcal{M}_i$.

Thus, we may think of the projection $\tilde{\pi}$ as a map from $\bar{V}_i \to \bar{S}_i$.

Let $\bar{K} > 0$ be such that $\mid \phi(z) \mid \leq \bar{K}$ for all $z \in \bar{V}_i$.

The function

$$h(z) = \begin{cases} \int_{\tilde{\pi}^{-1}(z)} \phi(w) d\nu_{\tilde{\pi}^{-1}(z)}(w) & \text{for } z \in S_i \\ 0 & \text{for } z \in \bar{S}_i \setminus S_i \end{cases}$$

is then bounded and measurable and its restriction to $S_i$ is continuous. Also, $\mid h(z) \mid \leq \bar{K}$ for all $z \in \bar{S}_i$.

Let $\bar{\mu}^i$ be the normalized restriction of $\bar{\mu}$ to $\bar{V}_i$.

We assert

$$\int_{\bar{S}_i} h(z) d(\tilde{\pi}_\star \bar{\mu}^i) = \int_{\bar{V}_i} \phi d\bar{\mu}^i \tag{102}$$

Since, $\bar{\mu}(\bar{\Sigma} \setminus \Sigma) = 0$, this tells us that the conditional measures of $\bar{\mu}$ with respect to $\mathcal{M}$ are the $\nu_C$ as required for Claim 2.

To prove (102), we let $\epsilon > 0$ be arbitrary, and we show

$$\mid \textstyle\int_{\bar{S}_i} h(z) d(\tilde{\pi}_\star \bar{\mu}^i) - \int_{\bar{V}_i} \phi d\bar{\mu}^i \mid \leq 5\epsilon \bar{K} \tag{103}$$

Let $\nu^i_{n_k}$ be the normalized restriction of $\nu_{n_k}$ to $\bar{V}_i$.

Since $\bar{V}_i$ is open and closed in $\bar{\Sigma}$, we have $\nu^i_{n_k} \to \bar{\mu}^i$ as $k \to \infty$.

Since $\tilde{\pi} : \bar{V}_i \to \bar{S}_i$ is continuous, we get $\tilde{\pi}_\star \nu^i_{n_k} \to \tilde{\pi}_\star \bar{\mu}^i$.

By (97), there is a compact subset $A_i \subset S_i$ such that, for large $k \geq 0$,

$$\nu^i_{n_k}(\bar{V}_i \setminus \tilde{\pi}^{-1}(A_i)) < \epsilon. \tag{104}$$

Since $h$ restricted to $A_i$ is continuous, we can use the Tietze extension theorem to find a continuous map $\bar{h} : \bar{S}_i \to \mathbf{R}$ such that $\mid \bar{h}(z) \mid \leq \bar{K}$ for all $z \in \bar{S}_i$, and $\bar{h}(z) = h(z)$ for $z \in A_i$.

Then,

$$\int_{\bar{S}_i} \bar{h} d(\tilde{\pi}_\star \nu^i_{n_k}) \to \int_{\bar{S}_i} \bar{h} d(\tilde{\pi}_\star \bar{\mu}^i).$$

By construction of $\bar{h}$, we then get, for large $k$,

$$\mid \textstyle\int_{A_i} h d(\tilde{\pi}_\star \nu^i_{n_k}) - \int_{A_i} h d(\tilde{\pi}_\star \bar{\mu}^i) \mid \leq 3\epsilon \bar{K}.$$

By (**),

$$\int_{A_i} h d(\tilde{\pi}_\star \nu^i_{n_k}) = \int_{\tilde{\pi}^{-1}(A_i)} \phi d(\nu^i_{n_k}),$$

The right side of this last equality differs from $\int_{\bar{V}_i} \phi d\nu^i_{n_k}$ by no more than $\epsilon \bar{K}$, and

$$| \int_{A_i} hd(\tilde{\pi}_\star \bar{\mu}^i) - \int_{S_i} hd(\tilde{\pi}_\star \bar{\mu}^i) | \le \epsilon \bar{K}.$$

Putting all these together gives (103) and completes the proof of Claim 2.

**Claim 3.** $\bar{\mu}$ is ergodic.

**Proof.**

This is a variant of the standard Hopf argument for geodesic flows in negatively curved Riemannian manifolds.

Let $\phi : \bar{\Sigma} \to \mathbf{R}$ be continuous. We show that $\bar{\mu}-$almost all forward time averages

$$\frac{1}{n} \sum_{k=0}^{n-1} \phi(\sigma^k z)$$

approach the same value.

Let

$$\phi_{for}(z) = \lim_{n \to \infty} \frac{1}{n} \sum_{k=0}^{n-1} \phi(\sigma^k z)$$

and

$$\phi_{bac}(z) = \lim_{n \to \infty} \frac{1}{n} \sum_{k=0}^{n-1} \phi(\sigma^{-k} z)$$

be the forward and backward limiting time averages of a point $z$.

From the Ergodic Theorem and standard arguments, there is a set $A_1 \subset \bar{\Sigma}$ of full $\bar{\mu}-$measure such that $z \in A_1$ implies $\phi_{for}(z), \phi_{bac}(z)$ exist and are equal. Also, since $\phi$ is continuous, $\phi_{for}$ is constant on stable sets and $\phi_{bac}$ is constant on unstable sets.

For each $z \in \Sigma$, let

$$W^s(z) = \bigcup_{n \ge 0} \sigma^{-n} W^s_{loc}(z)$$

be the global stable set of $z$.

Now, $\bar{\mu}-$almost any local unstable set $C$ is such that $\nu_C(A_1 \bigcap C) = 1$. Pick one such $C$ and let $\tilde{S}$ be the union of the global stable sets of points in $A_1 \bigcap C$. By the topological transitivity of the shift, the absolute continuity of the stable foliation $\mathcal{W}$ in $\tilde{Q}$, and the fact that the push forwards by $\pi$ of the conditional measures of $\bar{\mu}$ with respect to $\mathcal{M}$ are equivalent to the Riemannian measures on the local $F$-unstable manifolds, we get that $\nu_{C_1}(\tilde{S}) = 1$, for every local unstable set $C_1$. Hence, $\bar{\mu}(\tilde{S}) = 1$.

For any two points $z_1, z_2 \in \tilde{S}$, there are points $w_1, w_2 \in A_1 \bigcap C$ such that $z_1 \in W^s(w_1), z_2 \in W^s(w_2)$.

Then, $\phi_{for}(z_1) = \phi_{for}(w_1) = \phi_{bac}(w_1) = \phi_{bac}(w_2) = \phi_{for}(w_2) = \phi_{for}(z_2)$. This proves Claim 3.

**Claim 4:** $\lim_{n\to\infty} \nu_n = \bar{\mu}$.

**Proof.**

Let $\bar{\mu}_1$ be another subsequential limit of the sequence $\{\nu_n\}$. Substituting $\bar{\mu}_1$ for $\bar{\mu}$ in the preceding arguments gives that $\bar{\mu}_1$ is ergodic, shift invariant and $\bar{\mu}_1(\Sigma) = 1$.

Let $G_{\bar{\mu}}$ be the set of $\bar{\mu}-$generic points, and let $G_{\bar{\mu}_1}$ be the set of $\bar{\mu}_1-$generic points. Thus, for any continuous function $\phi : \bar{\Sigma} \to \mathbf{R}$,

$$\mathbf{a} \in G_{\bar{\mu}} \Rightarrow \frac{1}{n}\sum_{k=0}^{n-1} \phi(\sigma^k \mathbf{a}) \to \int \phi d\bar{\mu} \tag{105}$$

and

$$\mathbf{a} \in G_{\bar{\mu}_1} \Rightarrow \frac{1}{n}\sum_{k=0}^{n-1} \phi(\sigma^k \mathbf{a}) \to \int \phi d\bar{\mu}_1 \tag{106}$$

Ergodicity implies that $\bar{\mu}(G_{\bar{\mu}} \cap \Sigma) = 1 = \bar{\mu}_1(G_{\bar{\mu}_1} \cap \Sigma) = 1$.

If we show that

$$G_{\bar{\mu}} \cap G_{\bar{\mu}_1} \cap \Sigma \neq \emptyset \tag{107}$$

then, in view of (105) and (106), we get

$$\int \phi d\bar{\mu} = \int \phi d\bar{\mu}_1$$

for all continuous $\phi$, and Claim 4 follows.

For a given set $A \subset \Sigma$, let

$$W^s(A) = \bigcup_{\mathbf{a}\in A} W^s_{loc}(\mathbf{a})$$

We call $A$ *stably saturated* if $W^s(A) = A$. It is easy to see that both $G_{\bar{\mu}}$ and $G_{\bar{\mu}_1}$ are stably saturated.

The arguments in the proof of Claim 3 show that if $\bar{\mu}(A) = 1$ and $A \subset \Sigma$, then, for any local unstable set $C$, with Sinai measure $\nu_C$, we have $\nu_C(W^s(A)) = 1$. In particular,

$$\nu_C(G_{\bar{\mu}} \cap \Sigma) = \nu_C(W^s(G_{\bar{\mu}} \cap \Sigma)) = 1$$

Replacing $\bar{\mu}$ by $\bar{\mu}_1$ in the arguments of Claim 3 gives $\nu_C(G_{\bar{\mu}_1} \cap \Sigma) = 1$, as well. Thus, $\nu_C(G_{\bar{\mu}} \cap G_{\bar{\mu}_1} \cap \Sigma) = 1$ for any $C$ and (107) holds.

This completes the proof of Theorem 11.1.

**The construction of the SRB measure $\mu$.**

Let $\mu = \pi_\star \bar{\mu}$.

The measure $\mu$ is clearly an $F-$invariant and ergodic measure on $Q$.

There is a set $A \subset \tilde{Q}$ of full $\mu$ measure consisting of $\mu-$ generic points; i.e., $x \in A, \phi : Q \to \mathbf{R}$ continuous implies that $\frac{1}{n}\sum_{k=0}^{n-1} \phi(F^k x) \to \int \phi d\mu$.

Let $S$ be the union of the local stable manifolds of points $x \in A$. Clearly, each $x \in S$ is $\mu-$generic. We will show that $m(S) = 1$ (i.e. that $S$ has full Lebesgue measure in $Q$) to prove that $\mu$ is SRB.

Now, $\pi^{-1}(S)$ has full $\bar{\mu}-$measure in $\Sigma$. Hence, for some (in fact, $\bar{\mu}-$almost any) local unstable set $C \subset \Sigma$, we have $\nu_C(\pi^{-1}S) = 1$. This gives $\pi_\star\nu_C(S) = 1$. But $\pi_\star\nu_C$ is equal to the normalized Sinai measure on the local unstable manifold containing $S \bigcap \pi(C)$, and, hence, is equivalent to the Riemannian measure restricted to $S \bigcap \pi(C)$. This implies that $S \bigcap \pi(C)$ has full Riemannian measure in $\pi(C)$. Then, the absolute continuity of $\mathcal{W}$ gives $\rho_\gamma(S) = 1$ for every $K^u_\alpha$ curve $\gamma$, so Fubini's theorem gives $m(S) = 1$.

## 12 Further ergodic properties and an entropy formula

In this section we will study properties of the natural extension of the ergodic system $(F, \tilde{Q}, \mu)$. The first proposition identifies this natural extension with the system $(\sigma, \Sigma, \bar{\mu})$.

**Proposition 12.1** *The system $(\sigma, \Sigma, \bar{\mu})$ is isomorphic to the natural extension of the system $(F, \tilde{Q}, \mu)$.*

**Proof.** Since the map $F$ on $\tilde{Q}$ is not surjective, the meaning of this proposition is that there is a subset $\tilde{Q}_1$ of $\tilde{Q}$ of full $\mu$-measure such that $F(\tilde{Q}_1) = \tilde{Q}_1$, and the system $(\sigma, \Sigma, \bar{\mu})$ is isomorphic (mod 0) to the natural extension of the system $(F, \tilde{Q}_1, \mu)$.

Indeed, let $\tilde{Q}_1$ be the set of points $x \in \tilde{Q}$, such that there is a sequence $x_0, x_1, \ldots$ in $\tilde{Q}$ with $x_0 = x$ and $F(x_{n+1}) = x_n$ for all $n \geq 0$. It is easy to see that $F$ maps $\tilde{Q}_1$ onto itself. To see that $\mu(\tilde{Q}_1) = 1$, it suffices to show that $\bar{\mu}(\pi^{-1}\tilde{Q}_1) = 1$, and, since $\pi^{-1}\tilde{Q}$ has full $\bar{\mu}$ measure and $\bar{\mu}$ is $\sigma-$invariant, this follows from

$$\pi^{-1}(\tilde{Q}_1) \supset \bigcap_{n \geq 0} \sigma^n(\pi^{-1}\tilde{Q}) \tag{108}$$

To prove (108), let $\mathbf{a} \in \bigcap_{n \geq 0} \sigma^n(\pi^{-1}\tilde{Q})$, and let $x_0 = \pi(\mathbf{a}), x_n = \pi\sigma^{-n}\mathbf{a}$.

Since, $\sigma^{-n}\mathbf{a} \in \pi^{-1}\tilde{Q}$ for all $n \geq 0$, we have that $x_n = \pi\sigma^{-n}\mathbf{a} \in \tilde{Q}$ for each such $n$. On the other hand, $Fx_{n+1} = F\pi\sigma^{-n-1}\mathbf{a} = \pi\sigma\sigma^{-n-1}\mathbf{a} = \pi\sigma^{-n}\mathbf{a} = x_n$ for all $n \geq 0$. This shows that $x_0 = \pi\mathbf{a} \in \tilde{Q}_1$, so $\mathbf{a} \in \pi^{-1}\tilde{Q}_1$ which is (108). So, $\tilde{Q}_1$ is the required set.

The underlying set $\hat{Q}$ of the natural extension of $(F, \tilde{Q}_1, \mu)$ may be identified with the set of sequences $\bar{x} = (x_0, x_1, \ldots)$ in which each $x_n \in \tilde{Q}_1$ and $Fx_{n+1} = x_n$ for all $n \geq 0$.

Let $\xi = \{E_1, E_2, \ldots\}$ be the original collection of full height rectangles of $Q$. For any sequence $\bar{x} \in \hat{Q}$, the element $x_n$ is in the interior of a unique $E_{a_{-n}}$. Similarly, the point $F^n(x_0)$ is in the interior of a unique $E_{a_n}$. This enables us to define a map $\phi : \hat{Q} \to \Sigma$ by

$$\phi(\bar{x}) = \mathbf{a}$$

where, for each $n \geq 0$, $x_n \in int\ E_{a_{-n}}$ and $F^n x_0 \in int\ E_{a_n}$. Now, the verification that the map $\phi$ induces an isomorphism (mod 0) between the system $(\sigma, \Sigma, \bar{\mu})$ and the natural extension of $(F, \tilde{Q}_1, \mu)$ is straightforward, and we leave the details to the reader. QED.

Let $\zeta$ be the partition of $\Sigma$ into the sets $V_i$; i.e., the time 0 partition. Put $\eta = \bigvee_{i=-\infty}^{0} \sigma^i \zeta$. Then, the elements of $\eta$ coincide with the local stable sets $W^s_{loc}(\mathbf{a})$.

Moreover, we have that, mod zero, $\sigma\eta \succ \eta, \bigvee_n \sigma^n \eta$ is the point partition, and $\bigwedge_n \sigma^n \eta$ is the trivial partition $\{\Sigma\}$.

So, by definition, $(\sigma, \bar{\mu})$ is a K-system.

Then we state

**Proposition 12.2** *The map $(\sigma, \bar{\mu})$ is Bernoulli.*

We thank Dan Rudolph and Francois Ledrappier for useful conversations in connection with the proof of this proposition.

The following *Weak Markov* property was introduced in [11]. It was used to prove the Bernoulli property of Anosov flows ( see [4], [11]).

Let $\beta$ be any partition,

$$\beta_k^l = \bigvee_{k \leq i \leq l} \sigma^i \beta.$$

Given a collection of sets $P$, let us use $P^+$ for its union.

Say that $\beta$ is weak Markov (WM) if, for any $\epsilon > 0$, there is an integer $N = N(\epsilon)$, and collections $P = P(\epsilon)$ of atoms of $\beta_0^\infty$, $M = M(\epsilon)$ of atoms of $\beta_{-\infty}^0$ with the following properties.

1. $\bar{\mu}(P^+) > 1 - \epsilon$, and $\bar{\mu}(M^+) > 1 - \epsilon$.

2. For any $x_0^N \in \beta_0^N$,any $\bar{x}, \bar{y} \in P$ with $\bar{x} \bigcup \bar{y} \subset x_0^N$, and any subcollection $A$ of $M$ with $\bar{\mu}(A^+|\bar{y}) > 0$, one has

$$\mid \tfrac{\bar{\mu}(A^+|\bar{x})}{\bar{\mu}(A^+|\bar{y})} - 1 \mid < \epsilon \tag{109}$$

The proof of Proposition 2.2 in [11] shows that a finite weak Markov partition in a K-system is weakly Bernoulli in the sense of Friedman and Ornstein [5].

We will prove that the partition $\zeta$ is weak Markov. Then, arguments as in the proof of Proposition 2.2 in [11] give us that each of the finite partitions $\zeta_k = \{V_1, V_2, \ldots, V_k, \Sigma \setminus \bigcup_{i=1}^k V_i\}$ is also weakly Bernoulli. This implies that each factor map on $\Sigma / \bigvee_i \sigma^i \zeta_k$ is Bernoulli. Then, Theorem 5 in [9] gives that $(\sigma, \bar{\mu})$ is Bernoulli.

Thus, to prove Proposition 12.2, it suffices to show that the partition $\zeta$ of $\Sigma$ is weak Markov.

The corresponding $\zeta_0^\infty$ is the partition into local unstable sets $W^u_{loc}(\mathbf{a})$, and $\zeta^0_{-\infty}$ is the partition into local stable sets.

Given $\epsilon > 0$, let $n_0(\epsilon) >$ be large enough so that $\bar{\mu}(\bigcup_{\mid i \mid > n_0(\epsilon))} V_i) < \frac{\epsilon}{4}$.

For each $i$, let $z_i$ be a point in $V_i$, and let $A_i \subset W^s_{loc}(z_i), B_i \subset W^u_{loc}(z_i)$ be compact subsets so that the sets

$$D_i^u = \bigcup_{z \in A_i} W^u_{loc}(z), \; D_i^s = \bigcup_{w \in B_i} W^s_{loc}(w)$$

satisfy

$$\bar{\mu}(V_i \setminus D_i^u) < \frac{\epsilon}{2^{\mid i \mid + 2}}, \; \bar{\mu}(V_i \setminus D_i^s) < \frac{\epsilon}{2^{\mid i \mid + 2}}$$

Then, set $P = P(\epsilon) = \bigcup_{\mid i \mid \le n_0(\epsilon)} D_i^u, M = M(\epsilon) = \bigcup_{\mid i \mid \le n_0(\epsilon)} D_i^s$.

We have that $\bar{\mu}(P^+) > 1 - \epsilon, \bar{\mu}(M^+) > 1 - \epsilon$. Also, the set $Z_\epsilon = \bigcup_{\mid i \mid \le n_0(\epsilon)} D_i^s \bigcap D_i^u$ is compact.

Let $\mathbf{a}, \mathbf{b} \in V_i$, and let $W^u_{loc}(\mathbf{a}), W^u_{loc}(\mathbf{b})$ be their local unstable sets. Let $\pi_{\mathbf{a},\mathbf{b}}$ be the projection from $W^u_{loc}(\mathbf{a})$ to $W^u_{loc}(\mathbf{b})$ along the local stable sets in $V_i$. As $\mathbf{a}$ approaches $\mathbf{b}$ in $\Sigma$, the maps $\pi_{\mathbf{a},\mathbf{b}}$ approach the identity $\pi_{\mathbf{b},\mathbf{b}}$ and the measures $\rho_{\mathbf{a}}$ approach $\rho_{\mathbf{b}}$. Also, the densities $\bar{\xi}(\mathbf{a}, \mathbf{b})$ vary continuously with $\mathbf{a}, \mathbf{b}$ in $V_i$. On the compact set $Z_\epsilon$ the convergence and continuity above are uniform. Further, each $\bar{x} \in P$ is one of the sets $W^u_{loc}(\mathbf{a})$ and the conditional measure $\bar{\mu}(\cdot | \bar{x})$ is just the Sinai measure $\nu_{\bar{x}}$. If $\mathbf{a} \in \bar{x}, \mathbf{b} \in \bar{y}$ and $\bar{x} \bigcup \bar{y} \subset x_0^N \in \zeta_0^N$, then $a_j = b_j$ for $-N \le j \le 0$. For $N$ large, the measures $\bar{\mu}(\cdot|\bar{x}), \bar{\mu}(\cdot|\bar{y})$ have densities whose quotient is closer to 1 than $\epsilon$. These statements imply the Weak Markov property above. This completes the proof of Proposition 12.2.

**Entropy formula.**

It follows from our constructions that the measures of $V_i$ satisfy

$$c_1\delta_{i,min} < \bar{\mu}(V_i) < c_2\delta_{i,max} \tag{110}$$

for some positive constants $c_1, c_2$.

Since the partition $\zeta$ generates, we get

$$h_{\bar{\mu}}(\sigma) = \inf_n \frac{1}{n} H_{\bar{\mu}}(\zeta^0_{-n+1}) \leq H_{\bar{\mu}}(\zeta)$$

From condition G3 and (110), the last term is finite.

For $\mathbf{a} \in \Sigma$, and $n \geq 1$, let $V_{a_0 \dots a_{n-1}} = V_{a_0} \bigcap \sigma^{-1} V_{a_1} \bigcap \dots \bigcap \sigma^{-n+1} V_{a_{n-1}}$, and let $E_{a_0 \dots a_{n-1}}$ be the full height subpost of $Q$ defined in section 2.

Since $\sigma$ is ergodic with respect to $\bar{\mu}$, the Shannon-Breiman-Macmillan theorem gives a set $A$ with $\bar{\mu}(A) = 1$ such that $\mathbf{a} \in A$, implies

$$-\frac{1}{n} \log \bar{\mu}(V_{a_0 \dots a_{n-1}}) \to h_{\bar{\mu}}(\sigma) \tag{111}$$

Using that the conditional measures of $\bar{\mu}$ along local unstable sets have bounded densities relative to the measures $\rho_{\mathbf{a}}$, we see that there is a constant $K > 0$ such that, for $\mathbf{a} \in \Sigma$,

$$K^{-1} \min_{z \in W^s_{loc}(\pi \mathbf{a})} diam(\ell_z \bigcap E_{a_0 \dots a_{n-1}}) \leq \bar{\mu}(V_{a_0 \dots a_{n-1}}) \tag{112}$$

and

$$\bar{\mu}(V_{a_0 \dots a_{n-1}}) \leq K \max_{z \in W^s_{loc}(\pi \mathbf{a})} diam(\ell_z \bigcap E_{a_0 \dots a_{n-1}}). \tag{113}$$

This and Proposition 8.1 imply that, if $F^n(z) = (F^n_1(z), F^n_2(z))$, then there are a constant $K_1(\mathbf{a}) > 0$ and points $u_{n,1}(\mathbf{a}), u_{n,2}(\mathbf{a}) \in W^s_{loc}(\pi \mathbf{a})$ such that

$$\bar{\mu}(V_{a_0 \dots a_{n-1}}) | F^n_{1x}(u_{n,1}(\mathbf{a})) | \leq K_1(\mathbf{a}) \tag{114}$$

and

$$K_1(\mathbf{a})^{-1} \leq \bar{\mu}(V_{a_0 \dots a_{n-1}}) | F^n_{1x}(u_{n,2}(\mathbf{a})) | \tag{115}$$

By arguments like those in the proof of estimate (91), for $\bar{\mu}$ almost all $\mathbf{a}$, there is a constant $K_1(\mathbf{a}) > 0$ such that, for $z, w \in W^s_{loc}(\pi\mathbf{a}), n \geq 1$,

$$K_2(\mathbf{a})^{-1} \leq \frac{\mid F^n_{1x}(z) \mid}{\mid F^n_{1x}(w) \mid} \leq K_2(\mathbf{a}) \tag{116}$$

From (114), (115), (116) we get the existence of a constant $K_3(\mathbf{a}) >$, such that

$$K_3(\mathbf{a})^{-1} < \bar{\mu}(V_{a_0 \ldots a_{n-1}}) \mid F^n_{1x}(\pi\mathbf{a}) \mid < K_3(\mathbf{a}). \tag{117}$$

Thus there is a set $A$ with $\bar{\mu}(A) = 1$, so that if $\mathbf{a} \in A$, then

$$\lim_{n\to\infty} \frac{1}{n} \log \mid F^n_{1x}(\pi\mathbf{a}) \mid = h_{\bar{\mu}}(\sigma) \tag{118}$$

Since $\sigma$ is isomorphic to the natural extension of $F$, we have $h_\mu(F) = h_{\bar{\mu}}(\sigma)$.

Letting $A_1 = \pi(A)$, then, for $\mu$-almost all $z$ in $A_1$, we have

$$\lim_{n\to\infty} \frac{1}{n} \log \mid F^n_{1x}(z) \mid = h_\mu(F) \tag{119}$$

Taking $S$ to be the union of the stable manifolds of points in $A_1$, we get that $S$ has full Lebesgue measure in $Q$ and (119) holds for all $z \in S$.

But, for $z \in \tilde{Q}$, we have

$$\mid F^n_{1x}(z) \mid = \max(\mid F^n_{1x}(z) \mid, \mid F^n_{2x}(z) \mid) \leq \mid DF^n(z) \mid \leq (1+\alpha) \mid F^n_{1x}(z) \mid .$$

So, we have proved formula (6) and completed the proof of Theorem 3.1.

As a final remark, if $v = (v_1, v_2)$ is a unit vector in $K^u_\alpha$, then

$$\begin{aligned}
\left(1-\alpha^2\right) \mid F^n_{1x} \mid &\leq \mid F^n_{1x} v_1 + F^n_{1y} v_2 \mid \\
&= \mid DF^n(z) v \mid \\
&\leq \mid F^n_{1x} \mid \left(1+\alpha^2\right)
\end{aligned}$$

That is, for certain constants $C_1, C_2$, we have

$$C_1 \mid F^n_{1x} \mid \leq \mid DF^n(z)(v) \mid \leq C_2 \mid F^n_{1x} \mid$$

which, together with (119), implies formula (7).